\documentclass[11pt]{article}

\usepackage{fancyhdr}            
\fancypagestyle{plain}{
    \lhead{}
    \fancyhead[R]{\thepage}
    \fancyhead[L]{}
    
    \fancyfoot{}
}

\usepackage{array}
\usepackage[utf8]{inputenc} 
\usepackage[T1]{fontenc}    
\usepackage{xcolor}
\RequirePackage[colorlinks,citecolor=blue,urlcolor=blue]{hyperref}
\usepackage{nicefrac} 
\newcommand{\Spvek}[2][r]{%
  \gdef\@VORNE{1}
  \left(\hskip-\arraycolsep%
    \begin{array}{#1}\vekSp@lten{#2}\end{array}%
  \hskip-\arraycolsep\right)}

\def\vekSp@lten#1{\xvekSp@lten#1;vekL@stLine;}
\def\vekL@stLine{vekL@stLine}
\def\xvekSp@lten#1;{\def\temp{#1}%
  \ifx\temp\vekL@stLine
  \else
    \ifnum\@VORNE=1\gdef\@VORNE{0}
    \else\@arraycr\fi%
    #1%
    \expandafter\xvekSp@lten
  \fi}
\usepackage{algorithm}
\usepackage{algorithmic}
\usepackage{graphicx}
\usepackage{amsfonts}
\usepackage{amsmath}
\usepackage{url}
\usepackage{amsthm}
\usepackage{amssymb}
\usepackage{enumerate}
\usepackage{lmodern}
\usepackage{multirow}
\usepackage{subcaption}
\usepackage{graphicx}
\usepackage{graphicx}
\usepackage{amsfonts}
\usepackage{amsmath}
\usepackage[utf8]{inputenc}
\usepackage{url}
\usepackage{color}
\usepackage{amsthm}
\usepackage{amssymb}
\usepackage{enumerate}
 \usepackage{paralist}
\usepackage{graphicx}

\usepackage{bbm}
\usepackage{mathtools}
\RequirePackage{amsthm,amsmath,amsfonts,amssymb}
\RequirePackage{graphicx}
\usepackage{microtype}  
\usepackage{times}
\usepackage{mathrsfs}
\usepackage{mathtools}
\usepackage{transparent}
\usepackage{graphicx}
\usepackage{wrapfig}
\usepackage{relsize,exscale}
\usepackage{nameref}

\usepackage{url}
\usepackage{amsthm}
\usepackage{amssymb}
\newtheorem {Proposition}{Proposition}[section]
\newtheorem {Lemma}[Proposition] {Lemma}
\newtheorem {Theorem}[Proposition]{Theorem}

\newtheorem*{Theorem*}{Theorem}

\newtheorem {Remark}[Proposition]{Remark}

\newtheorem {Definition}{Definition}[section]

\def\N{\mathbb{N}}

\def\R{\mathbb{R}}
 
\def\E{\mathbb{E}}
\def\P{\mathbb{P}}

\usepackage{relsize}
\usepackage{tikz}

\makeatletter
\def\namedlabel#1#2{\begingroup
   \def\@currentlabel{#2}%
   \label{#1}\endgroup
}
\makeatother

\def\P{\mathbb P}
\def\E{\mathbb E}

\def\N{\mathbb N}
\def\H{\mathcal{H}}

\def\P{\mathcal{P}}
\def\Pac{\P^{\rm {\,}a.c.}}
\def\Pdc{\P^{\rm {\,}d.c.}}
\def\Pl{\P^{\,\ell}}
\def\PH{\P^{\rm H}}

\usepackage{minitoc}
\usepackage{enumerate}
\usepackage{bbm}
\usepackage[T1]{fontenc}    
\usepackage{lmodern}        
\usepackage{hyperref}       
\usepackage{url}            
\usepackage{booktabs}       

\usepackage{amsfonts}       

\usepackage{nicefrac}       
\usepackage{microtype}      
\usepackage{lipsum}
\usepackage{graphicx}
\usepackage{color}
\graphicspath{ {./images/} }
{
     \theoremstyle{plain}
     
}
\usepackage[sectionbib,sort,compress]{natbib}
\usepackage{chapterbib}

\begin{document}

\begin{center}
{\Large Monotone Measure-Preserving Maps in Hilbert Spaces:\\ Existence, Uniqueness, and Stability} \\
\vspace{0.5cm} 
Alberto Gonz\'alez-Sanz,$^{1}$ 
Marc Hallin,$^{2}$ and 
Bodhisattva Sen$^{3}$
\\
\vspace{0.2cm} 
$^1$ Institut de Mathématiques de Toulouse,  Université Paul Sabatier, Toulouse, France\\
$^2$ ECARES and Département de Mathématique Université libre de Bruxelles, Brussels, Belgium\\
$^3$ Department of Statistics, Columbia University, New York, United States 
\end{center}

\begin{abstract}
The contribution of this work is twofold. The first part deals with a Hilbert-space version of McCann's celebrated result on the existence and uniqueness of monotone measure-preserving  maps: given two probability measures $\rm P$ and $\rm Q$ on a separable Hilbert space $\H$ 
where $\rm P$ does not give mass to ``small sets''  (namely, {\it Lipschitz hypersurfaces}), we show, without imposing any moment assumptions, that there exists a gradient of  convex function $\nabla\psi$  pushing  ${\rm P} $ forward to~${\rm Q}$.   In case $\H$ is infinite-dimensional, ${\rm P}$-a.s.~uniqueness  is not guaranteed, though. 
If, however,~${\rm Q}$ is  boundedly supported (a natural assumption in several statistical applications), then this  gradient is~${\rm P}$-a.s.~unique. In the second part of the paper, we establish    stability results for transport maps in the sense of uniform convergence over compact ``regularity sets''.   As a consequence, we obtain a central limit theorem for the fluctuations of the optimal quadratic transport cost in a separable Hilbert space. 
\end{abstract}


\section{Introduction}\label{intro}
\subsection{Brenier and McCann}
Two seminal results  had a major impact on the recent surge of interest in measure transportation methods and their applications. The first one is the polar factorization theorem~\citep{brenier1991polar},  
associated with the name of Yann Brenier, although several authors~\citep{Cuesta1989NotesOT, RR90} independently contributed  partial versions of the same result. The second one~\citep{McCann}, which  extends the generality of Brenier's theorem by relaxing the moment conditions, is due to Robert McCann. 

Let ${\rm P}$ and ${\rm Q}$  belong to the family $\P(\R^d)$ of Borel probability measures over $\R^d$, for $d \ge 1$. Under its most usual version  (see, e.g.,  {Theorem~2.12 in  \cite{villani2003topics}}),  McCann's theorem states that, 
 for ${\rm P}$ in the Lebesgue-absolutely-continuous subfamily $\Pac(\R^d) \subset \P(\R^d)$,  there exists a~$\rm P$-a.s.\ unique gradient of  convex function~$\nabla\psi$~pushing  ${\rm P}$ forward to ${\rm Q}$ (notation:~$\nabla\psi \# {\rm P}={\rm Q}$); in case~${\rm P}$ and ${\rm Q}$ admit finite moments of order two, that gradient, moreover, is the {$\rm P$-a.s.\ unique} solution  of the quadratic optimal transport problem 
\begin{equation}
    \label{OTIntro}
\mathcal{T}_2({\rm P},{\rm Q})\coloneqq\inf_{T\#{\rm P}={\rm Q}} \int \|T(x)-x \|^2\,d {\rm P} (x).
\end{equation} 
 
Actually, \cite{McCann} established this result under the weaker assumption that $\rm P$ belongs to the class $\PH(\R^d) \supsetneq \Pac(\R^d)$ of probability measures vanishing on all Borel sets  with Hausdorff dimension $(d-1)$. McCann's result constitutes a substantial extension of Brenier's theorem which, under the restrictive assumption of finite second-order moments,\footnote{Finite second-order moments and ${\rm P} \in {\mathcal P}^{\text{{\,}a.c.}}(\R^d)$ are sufficient (see {Chapter 2 in\cite{villani2003topics})} for Brenier's result.   
\cite{brenier1991polar}, however, had further additional assumptions involving, e.g., the density and the support of $\rm P$, which are not necessary.} only implies that a gradient of convex function is the $\rm P$-a.s.\ unique solution of the  transport problem \eqref{OTIntro}. In his proofs, McCann adopted geometric ideas rather than analytical ones to prove his result; as commented in \cite{GaMc}, his argument can be related to that of Alexandrov's uniqueness proof for convex surfaces with prescribed Gaussian curvature.  

\subsection{Measure transportation in Hilbert spaces}
Now suppose that $\H$ is a separable Hilbert space with inner product $\langle \cdot, \cdot \rangle$ and induced norm $\|\cdot \|$. A very natural question is 
\textit{``Can we extend McCann's theorem~\citep{McCann} from the finite-dimensional real space $\R^d$ to the case of a general separable {\it Hilbert space} $\H$?''} In other words, given two probability measures $\rm P$ and $\rm Q$ in the family $\P(\H)$  of all Borel probability measures on~$\H$ such that $\rm P$ does not give mass to ``small sets'',  does there exist a unique gradient of convex function~$\nabla\psi$  pushing  ${\rm P} $ forward to ${\rm Q}$? 
\vspace{0.05in}

Theorem~\ref{Theorem:Maccan} in Section~\ref{sec:Exist-Unique} provides an affirmative answer to the above question by showing that, for any separable Hilbert space $\H$, provided that $\rm P$ gives  zero mass to so-called \textit{Lipschitz surfaces}, there exists a convex function $\psi$  the gradient $\nabla\psi$ of which  pushes ${\rm P}$ forward to ${\rm Q}$. Under the additional assumption that the support of ${\rm Q}$ is bounded, we further show that such a gradient of convex function $\nabla \psi$ is ${\rm P}$-a.s.~unique. 

To the best of our knowledge, the first results on the existence of  optimal transport  mappings  in Hilbert spaces are due to \cite{Cuesta1989NotesOT} who 
 prove the existence of solutions of~\eqref{OTIntro} (for ${\rm P}$ and ${\rm Q}\in \P(\H)$) under the following assumption on $\rm P$: for any basis~$\{ e_i\}_{i\in \N} $ of $\H$ and for any set $E \subset \H$ with ${\rm P}(E)>0$, there exists $a\in \H$ such that~$\mu_1(\{t\in~\!\R:~a+~\!t\, e_i~\!\in~\!E\})>~\!0$ for all $i\in \N$, where $\mu_1$ denotes the univariate Lebesgue measure. A probability measure satisfying this assumption, in particular, gives no mass to Aronszajn null sets.\footnote{Recall that { $ E\subset {\H}$} is an Aronszajn null set (cf.~\cite{Csrnyei1999AronszajnNA})  if there exists a complete sequence $\{ e_i\}_{i\in \N}\subset \H$ such that~$ E$  can be written as a union of Borel sets $\{E_i\}_{i \ge 1}$ such that each $E_i$ is null on every line in the direction $e_i$, i.e., for every $a\in \H$, $\mu_1(\{t\in \R: \ a+t\, e_i \in E_i\})=0$ for all $
 i\in \N$.} A uniqueness result for the same problem is established in \citet[Theorem~6.2.10]{AmbrosioGradient} under the additional assumption of finite second-order moments  for ${\rm P}$ and ${\rm Q}$. The argument for that uniqueness result takes advantage of the strict convexity of the functional in the right-hand side of~\eqref{OTIntro}, and is therefore  helpless in the absence of finite second-order moments. Thus, so far, no McCann extension of Brenier-type results is available in the general Hilbert space setting. 

Let us now comment on the main hurdles encountered in proving Theorem~\ref{Theorem:Maccan}.~\citet{McCann} showed the existence of such a $\nabla \psi$ pushing forward ${\rm P}$ to ${\rm Q}$, when $\H = \R^d$, by using a Rademacher-type result (see~\citet{Anderson1952}) which implies that a lower semi-continuous (l.s.c.)~convex function~$\varphi:\H \to (-\infty, \infty]$ is continuous on the interior of its domain and differentiable except on a set of Hausdorff dimension $d-1$ in $\operatorname{dom}(\varphi)$.\footnote{Here $\operatorname{dom}(\varphi) := \{x \in \H: \,\varphi(x) \in \R\}$ denotes the domain of $\varphi$.} Although there are infinite-dimensional extensions of the above result (see e.g.,~\cite{Zajk1979} or~\citet[Theorem 6.2.3]{AmbrosioGradient}), these results assume continuity and/or a local Lipschitz property of the underlying l.s.c.~convex function $\varphi$. Now, when $\H$ is infinite dimensional, there exists proper l.s.c.~convex functions $f:\H \to (-\infty, \infty]$ discontinuous at every point of $\H$ such that $\nabla f$  pushes forward a non-degenerate Gaussian distribution to another; see Remark~\ref{remark:nonGauss} for the details. We circumvent this difficulty by showing both existence and uniqueness of such a $\nabla \psi$ pushing forward ${\rm P}$ to a boundedly supported target measure, then creating a sequence of distributions with increasing but bounded supports to approximate ${\rm Q}$. Note that when the target measure is boundedly supported, following the arguments in~\citet[p.~147]{AmbrosioGradient}, we can show that $\nabla \psi$ exists with $\psi$  agreeing  ${\rm P}$-a.s.\ with a continuous convex function~$\bar{\psi}$. As a consequence, we can assume that $\psi$ is continuous in $\H$ when ${\rm Q}$ is boundedly supported.

To prove the uniqueness of $\nabla \psi$ under the assumption that $\operatorname{supp}({\rm Q})$ is bounded, we show that if two continuous convex functions~$f$ and $g$ have different gradients at a point $x \in \H$ (that is,~$\nabla f(x) \ne \nabla g(x)$), then there exists a neighborhood $\mathcal{U}_x$ of $x$ such that $\mathcal{U}_x \cap \{f=g\}$ belongs to the class of  {\it Lipschitz hypersurfaces} which, under the assumption that ${\rm P} $ does not give mass to such a class of sets, i.e., ${\rm P} \in \Pl(\H)$ (see Definition~\ref{defn:Lipschitz-Surfaces}-$(ii)$), is  a $\rm P$-null set. As a consequence, if $x \in \operatorname{supp}({\rm P})$, the set~$\mathcal{V}_x \coloneqq \mathcal{U}_x \cap \{f \ne g\}$ has strictly positive $\rm P$-measure. A contradiction is now obtained by noting that ${\rm P}(\nabla f \in \partial g(\mathcal{V}_x)) \ne {\rm P}(\nabla g \in  \partial g(\mathcal{V}_x))$, which makes $\nabla f \# {\rm P} = \nabla g \# {\rm P} = {\rm Q}$ impossible. 

Note that, in particular, for $\H=\R^d$, $\Pl(\R^d)\supseteq \PH(\R^d)$, and, for general $\H$,  any non-degenerate Gaussian measure belongs to $\Pl(\H)$ (see Section~\ref{sect:nullsets}).

\subsection{Stability of Hilbertian transport maps}
The second objective of this paper (Section~\ref{sec:stability}) is a characterization of the stability properties of the transport map $\nabla\psi$---a problem that has not been considered so far in   infinite-dimensional spaces.

The most general results in the  finite-dimensional case are due to  \cite{GhosalSenAOS},~\cite{Barrio2022NonparametricMC}, and \cite{Segers2022GraphicalAU}. Being based on the Fell topology, which does not have nice properties in non-locally compact spaces, the techniques used by these authors do not extend to general Hilbert spaces. Let us briefly describe the stability result when $\H = \R^d$. Let $\{{\rm P}_n\}_{n\in \N}$ and~$\{{\rm Q}_n\}_{n\in \N}$  be two sequences of probability measures on~$\R^d$ such that~${\rm P}_n\xrightarrow{w}{\rm P}$ and~${\rm Q}_n\xrightarrow{w}~\!{\rm Q}$, as~$n\to\infty$, where $\xrightarrow{w}$ denotes weak convergence of probability measures. Recall that the subdifferential of a l.s.c.~convex function~$\psi:\H \to (-\infty, +\infty]$ is defined  as
$$\partial \psi \coloneqq \{(x,y) \in \H \times \H: \psi(x) + \langle y, z - x \rangle \le \psi(z)\; \mbox{ for all }\; z \in \H\}.$$
Denote by $\Pi({\rm P}_n,{\rm Q}_n)$ the family of distributions in ${\mathcal P}(\H\times \H)$ with marginals~${\rm P}_n$ and ${\rm Q}_n$, and let $ \gamma_n\in \Pi({\rm P}_n,{\rm Q}_n)$ be such that~$\operatorname{supp}(\gamma_n)\subseteq  \partial \psi_n$ for some l.s.c.~convex function~$\psi_n:\H \to (-\infty, +\infty]$. Further, let $ \psi: \R^d \to (-\infty, +\infty]$ be a proper l.s.c.\ convex function such that~$ \nabla \psi$ pushes $\rm P$ forward to~$\rm Q$. 
Then,  for any compact subset~$K$ of~$\operatorname{dom}(\nabla \psi)\cap \operatorname{int}(\operatorname{supp}({\rm P}))$,\footnote{For notational convenience, we write $\sup_{(x,y)\in \partial \psi_n, x\in K } \| y-\nabla\psi(x)\| \coloneqq \sup_{x \in K} \sup_{y \in \partial \psi_n(x)} \| y-\nabla\psi(x)\|$.} (here $\operatorname{int}(\cdot)$,  $\operatorname{dom}(\cdot)$, and $\operatorname{supp}(\cdot)$ stand for the interior of a set, the domain of a function, and the support of a distribution, respectively)   
\begin{equation}\label{stabfin}
    \sup_{(x,y)\in \partial \psi_n, x\in K } \| y-\nabla\psi(x)\|  \longrightarrow 0\quad\text{as $n\to\infty$}.
\end{equation}
In this Euclidean setting, \eqref{stabfin} holds  without any assumption on ${\rm Q}$. 

Section \ref{sec:stability} extends this finite-dimensional stability result to arbitrary {separable} Hilbert spaces. A Hilbert space $\H$, however, has two useful topologies: the {\it strong topology} under which $x_n\to x $ if and only if $\| x_n-x \|\to 0$ and the {\it weak topology} under which $ x_n\rightharpoonup x$ if and only if $ \langle h, x_n\rangle\to \langle h, x\rangle$ for all $h\in \H$. In the finite-dimensional case, these two topologies coincide, but they are distinct in the infinite-dimensional case.   Due to the fact that the map $\nabla \psi$ is only a.s.\ {\it strong-to-weak continuous}---it is mapping strongly convergent sequences to weakly convergent ones---in the set of differentiability points of $\psi$ (see~\citet[Theorem 21.22]{bauschkeMonotoneHilbert} and Section~\ref{sec:stability} for  formal definitions), 
we cannot expect  convergence in norm  as in \eqref{stabfin} to hold in general $\H$: our Theorem~\ref{Theorem stability} yields,  
for any strongly\footnote{By {\it strongly compact}  we mean compact with respect to the strong (norm) topology.} compact set~$K\subseteq  \operatorname{dom}(\nabla \psi)\cap \operatorname{int}(\operatorname{supp}({\rm P}))$,
\begin{equation}
    \label{StabilityIntro}
    \sup_{(x,y)\in \partial \psi_n, x\in K } \langle y-\nabla\psi(x), h \rangle  \longrightarrow 0\quad\text{as $n\to\infty$},
\end{equation}
for any $h\in \H$, i.e., stability in the weak topology. In the infinite-dimensional case we show via an example (see remark~(a) in Section~\ref{sec31})  that, without the assumption that $\operatorname{supp}({\rm Q})$ is bounded,~\eqref{StabilityIntro} can fail.

Our proof strategy for Theorem~\ref{Theorem stability} is as follows. 
We first prove (Lemma \ref{lem:Coupling}) the stability of the optimal (cyclically monotone) couplings $\gamma_n \in \Pi({\rm P}_n, {\rm Q}_n)$.  As a second step, we establish the convergence of the subdifferential $\partial \psi_n$ as a set-valued map; since we are dealing with cyclically monotone set-valued maps, {\it graphical convergence} in the sense of Painlev\'{e}-Kuratowski \citep[p.~111]{rockafellar2009variational} provides the appropriate framework. Mr\'{o}wka's theorem (Lemma~\ref{Mrowka}) then guarantees the existence of a graphical limit along subsequences. We show (Lemma~\ref{Lemma:LimitCyclically}) that the cyclical monotonicity of $\partial \psi_n$ is preserved in this graphical limit. The final step 
establishes that this limit, moreover, is contained in $\partial \psi$. This is achieved with Lemma~\ref{Lemma:uniqueness},  of independent interest, where we show  that  if the subdifferentials of two  convex functions  coincide on a dense subset of some  convex open set $\mathcal{B} \subset \H$, then they coincide on the entire set $\mathcal{B}$. 

Theorem~\ref{Theorem stability} also entails the stability of the potentials (whenever they are unique, up to additive constants) defining the transport maps. The proof   follows along similar lines as in the Euclidean case but is more involved  due to the fact that $\H$ may not be locally compact and hence Arzel\'{a}-Ascoli (see~\citet[Theorem 4.25]{Brezis}) may not apply. 
To overcome this,  we take advantage of the fact that, since ${\rm P}$ is tight, we can restrict the study of the convergence of $\psi_n$ to compact sets with arbitrarily large ${\rm P}$-probability. 
 
 Finally, denoting by ${\rm P}_n\coloneqq \frac{1}{n} \sum_{i=1}^n \delta_{X_i}$  the empirical distribution of a random sample $X_1,\ldots, X_n$ from ${\rm P}$, we obtain, in Theorem~\ref{TCL}, the central limit result 
$$\sqrt{n}\left(\mathcal{T}_2({\rm P}_n,{\rm Q})-\mathbb{E}\mathcal{T}_2({\rm P}_n,{\rm Q})\right)\stackrel{w}{\longrightarrow} N(0, \sigma^2_2({\rm P},{\rm Q})) $$
for the fluctuations of the squared $2$-Wasserstein distance {$\mathcal{T}_2({\rm P}_n,\rm Q)$} about its mean.\footnote{Here we assume that $\rm P\in \Pl({\H})$ admits finite fourth-order moments and ${\rm Q} \in \P({\H})$ has bounded support.} This result extends to general Hilbert spaces the finite-dimensional result by \cite{delBarrioLoubes19}.

\subsection{Statistical applications: Hilbert space-valued ``center-outward'' distribution and rank functions}

Observations, in a variety of statistical and machine learning problems, increasingly often take values in more complex spaces than $\R^d$ and infinite-dimensional Hilbert-space-valued observations (\cite{Small-McLeish-Book-1994}) nowadays are frequent---in functional data analysis (\cite{Horvath-Kokoszka-Book-2012, Hsing-Eubank-2015,Kokoszka-Reimherr-2017}), in the so-called {\it kernel methods} for general pattern analysis (e.g., in object-oriented data analysis, see~\cite{Marron-Alonso-2014}), in kriging theory for random fields~(\cite{Menafoglio-JMA-2016}), in shape analysis~(\cite{jayasumana2013}), etc. Moreover, the use of measure-transportation-based techniques to analyze such complex data  has also become increasingly important, with direct implications in several problems involving  Hilbert-space-valued data, such as  two-sample testing~\citep{Cuesta-Albertos-2006},  independence testing~\citep{Lai-Et-Al-2021},  quantile estimation~\citep{spatial-AoS-2014},  data depth~\citep{On-data-depth-2014}, etc. 


The finite second-order moment assumption required, e.g., in~\citet[Theorem~6.2.10]{AmbrosioGradient},  needs not to be satisfied, though.
This makes 
the McCann-type generalization in Theorem~\ref{Theorem:Maccan} essential and important in  statistical problems with Hilbert-space-valued observations.

A major statistical application of measure transportation in the $d$-dimensional Euclidean space is the definition of  multivariate concepts of ``center-outward'' distribution, rank and quantile functions and their empirical counterparts satisfying all the properties that make their univariate counterparts fundamental tools for statistical inference. These, in particular, allow for the construction of  {\it rank}-based methods (distribution-free rank-based testing and R-estimation) in $\R^d$. 

Recall that the distribution function of a continuous univariate random variable $X\sim\rm P$ is defined as~$x\mapsto F(x)\coloneqq {\rm} \mathbb{P}(X\leq x)$ for $x \in \R$. This distribution function $F$ actually  is  the unique gradient of a convex function pushing $\rm P$ forward to the uniform distribution over $(0,1)$, which exists and is $\rm P$-a.s.\ unique irrespective of the existence of any moments. Similarly, given a sample~$X_1,\ldots, X_n \sim {\rm P}$, the empirical distribution function is the transport map pushing the empirical measure of the $X_i$'s forward  to $\frac{1}{n} \sum_{i=1}^n \delta_{i/(n+1)}$---a natural discretization of the uniform distribution over $(0,1)$---while minimizing the quadratic transportation cost~\eqref{OTIntro}. 
This measure-transportation-based characterization has been used successfully to define multivariate versions of the concepts of ``center-outward'' distribution, multivariate rank and quantile functions, with the Lebesgue uniform over the unit cube \citep{chernoetal17, BodhisattvaDeb} or the spherical uniform over the unit ball \citep{Hallin2020DistributionAQ, FIGALLI2018413, DELBARRIO2020104671,Barrio2023RegularityOC} playing the role of the  reference distributions $\rm Q$;\footnote{Due to its strong symmetry properties,  the spherical uniform over the unit ball, unlike the Lebesgue uniform over the unit cube,  induces  adequate notions of quantile function and quantile regions.} these distributions  are boundedly supported, hence enter the realm of our uniqueness results. 

The multivariate rank (or ``center-outward'' distribution) function $\mathbf{F}$ of 
$ {\rm P}\in \P^{\text{a.s.}}(\R^d)$ is then defined as the unique gradient of a convex function pushing $ {\rm P}$ forward to the  reference distribution~$ {\rm Q}$. The essential properties of $\mathbf{F}$,  matching the properties of the traditional univariate distribution function, are: (i) distribution-freeness (as $\mathbf{F}(X) \sim {\rm Q}$ if $X \sim {\rm P}$); (ii) $\mathbf{F}$ entirely characterizes $ {\rm P}$; (iii) the (natural) empirical version of~$\mathbf{F}$ is uniformly consistent at its continuity points (see \cite{Hallin2020DistributionAQ} for details). It is worth mentioning that in this context, a  sensible definition of the concept of a ``center-outward'' distribution or rank function cannot be subjected to the existence of finite moments;   a McCann-type approach, thus, as opposed to  Brenier's  finite-second moment one,  is essential.

These new concepts have been applied to a wide range of inference problems such as vector independence and goodness-of-fit testing \citep{Shi2022, ShiSpear, BodhisattvaDeb, GhosalSenAOS}, testing for multivariate symmetry~\citep{HuangSen2023},  distribution-free rank-based testing and R-estimation for VARMA models \citep{HLLJASA, HLLBern, HLPortm}, multiple-output linear models and MANOVA \citep{HHH},  multiple-output quantile regression  \citep{Barrio2022NonparametricMC}, definition of multivariate Lorenz functions \citep{HMLorenz}, etc.; see \cite{AnnR} for a recent survey.  

Defining adequate concepts of a ``center-outward'' distribution or multivariate rank function in dimension $d>1$ has been an open problem in the statistical literature for about half a century. Many definitions have been proposed in this direction, including the many notions  of {\it statistical depth} following Tukey's celebrated concept (see~\cite{TUKEY}). None of these definitions, however,  yield the  essential  properties~(i)-(iii) 
of the traditional univariate concept mentioned above.  By establishing the existence and uniqueness of the gradient of a convex function pushing $ {\rm P} \in \Pl(\H)$, where $\H$ is a separable Hilbert space, forward to a boundedly supported ${\rm Q}$, our   Theorem~\ref{Theorem:Maccan}, which does not require $\rm P$ to admit  any moments, is a first step in  the direction of extending this measure-transportation-based approach to multivariate rank  (or ``center-outward'' distribution)  function from Euclidean spaces to general separable Hilbert spaces.


\section{Existence and uniqueness of monotone measure-preserving maps in Hilbert spaces}\label{Section:McCan}
\subsection{Preliminaries: definitions and notation}
\subsubsection{Some results from convex analysis}
Throughout, denote by  $\H$  a separable Hilbert space with inner product $\langle \cdot, \cdot \rangle$ and norm   $\|\cdot \|$. Two topologies can be considered for $\H$: the {\it strong topology}, under which $x_n\to x $ as $n \to \infty$ (where~$\{x_n\}_{n \ge 1} \subset \H$) if and only if $\| x_n-x \|\to 0$, and the {\it weak} one,  under which $ x_n\rightharpoonup x$ if and only if $ \langle h, x_n\rangle\to \langle h, x\rangle$ for all $h\in \H$.  The weak and strong topologies in $\H$ generate the same Borel $\sigma$-algebra (see \cite{Edgar1977}). 

Recall that a set $\Gamma\subseteq \H\times \H$ is said to be \textit{cyclically monotone} if, for all $n\in \N$ and all~$\{({x}_k,{y}_k)\}_{k=1}^{n} \subseteq \Gamma $, letting   $y_{n+1}=y_1$,
\begin{align}\label{Cyclically}
	\sum_{k=1}^n \langle x_{k},y_{k+1}-y_k\rangle\leq 0.
\end{align}
A Borel probability measure $\gamma \in \P(\H\times \H)$ is said to have ${\rm P} \in \P(\H)$ and ${\rm Q} \in \P(\H)$ as its (left
and right, respectively) \textit{marginals} if~$\gamma (A\times \H)={\rm P}(A) $ and $\gamma (\H\times B)={\rm Q}(B) $ for all Borel sets~$A,B\subseteq\H$. The family  of  $\gamma$'s having marginals ${\rm P}$ and ${\rm Q}$ is denoted by $\Pi({\rm P}, {\rm Q})$.  

Cyclically monotone sets and convex functions are related  in the following sense. Let  $\Gamma \subseteq \H\times \H$ be cyclically monotone. Theorem~B in~\cite{RockafellarMaximalMonot} establishes the existence of a proper lower semi-continuous (l.s.c.)  convex function $f:\H\to (-\infty, +\infty]$ such that $\Gamma$ is contained in the subdifferential $ \partial f$ of $f$: 
$$\Gamma\subseteq \partial f\coloneqq \{(x,y)\in \H\times \H: f(x)-f(x')\leq  \langle y, x-x'\rangle, \ \text{for all} \ x'\in \H \}. $$
 Without any loss of generality, the subdifferential $ \partial f$ can be assumed to be maximal monotone in the sense that $\partial f\subseteq\partial g $ for some other  proper l.s.c.\   convex function $g$  implies $\partial f=\partial g $. 

Slightly abusing notation, for each $x\in \H$, call $\partial f(x)\coloneqq \{y\in \H: \ (x,y)\in \partial f\}$ the  \emph{subdifferential at $x$} of $f$. The mapping $x\mapsto \partial f(x)$ is generally multi-valued. For a set $A\subseteq \H$, write 
\[ \partial f (A)\coloneqq \bigcup_{a\in A}\partial f(a).
\]
In case $\partial f(x) $ is a singleton, denote by $\nabla f(x)$ its unique element. 

When the Hilbert space $\H$ is not finite-dimensional, some of the familiar properties of convex functions no longer hold. For instance, the continuity of a convex $f$ in its domain is no longer guaranteed: 
$$ \operatorname{dom}(f)\coloneqq\{ x\in \H: \ f(x)\in \R\}\;\;\neq \;\; \operatorname{cont}(f)\coloneqq\{ x\in \H: \ \text{$x\mapsto f(x)$ is continuous}\}.$$
However, when $f$ is a  proper l.s.c.\   convex function, $\operatorname{int}(\operatorname{dom}(f))$ and $\operatorname{cont}(f)$ coincide ~\citep[Corollary 8.30]{bauschkeMonotoneHilbert}. Moreover, in that case, Proposition 16.14 (Ibidem) 
yields 
\begin{equation}
    \label{eq:contcontdom}
   \operatorname{int}(\operatorname{dom}(f))= \operatorname{cont}(f) \subseteq \operatorname{dom}(\partial f)\coloneqq \{ x\in \H: \ \partial f(x)\neq \emptyset \} \subseteq \operatorname{dom}(f),
\end{equation}
provided that $  \operatorname{int}(\operatorname{dom}(f))\neq \emptyset $.
Note that the domain of differentiability of a convex function $f$, denoted by
\begin{equation}\label{eq:Dom-Grad-f}
\operatorname{dom}(\nabla f)\coloneqq \{ h\in \H: \partial f (h) \; \text{is  a singleton} \}, 
\end{equation}
in infinite-dimensions, differs from\footnote{Recall from \cite{Shapiro1990OnCO} that a proper function $f:\H\to (-\infty, +\infty]$ is {\it Fréchet-differentiable} at~$h_0\in\H$ if there exists $a^*\in\H $ such that 
  $\lim_{h\to 0 } {\|f(h_0+h)-f(h_0)-\langle a^*, h\rangle \|}/{\|h\|}=0$. } 
$$ \operatorname{dom}_{\text{\rm Fr}}(\nabla f )\coloneqq \{ h\in \H:  f \; \text{is Fréchet-differentiable at $h$} \}. $$

The following lemma gives some basic continuity properties of the subdifferential of a proper l.s.c.~convex function defined on a Hilbert space---the continuity of the subdifferential $\partial f$ of $f$ depends on the kind of differentiability considered. Part~$(i)$ of Lemma~\ref{Lemma:strongtoweak} is a direct consequence  of the fact that, in the product space~$\H\times \H$ with  the first $\H$ factor equipped with the weak topology and the second one with the strong topology, the subdifferential $\partial f$ is a closed locally bounded\footnote{That is, for any $x\in \operatorname{int}(\operatorname{dom}(f))$, there exists a ball $\mathcal{B}(x, \epsilon)$ centered at $x$ such that $\partial f(\mathcal{B}(x, \epsilon))$ is bounded.}  set   (see Propositions 16.26 and 16.14 in \cite{bauschkeMonotoneHilbert}). Parts $(ii)$ and $(iii)$ can be found in Propositions 17.32  and  17.33 (Ibid.).  
\begin{Lemma}\label{Lemma:strongtoweak}
    Let $ f :\H\to (-\infty, +\infty]$ be a  proper l.s.c.\   convex function,  $x\in \operatorname{int}(\operatorname{dom}(f))$, and~$\{x_n\}_{n\in \N} \subset \H$ be a sequence such that $x_n\to x$ as $n \to \infty$. Then, 
    \begin{enumerate}
    \item[(i)]for any sequence $\{y_n \}_{n\in \N}$ with $y_n\in \partial f (x_n)$, there exists a subsequence weakly converging to~$y\in \partial f(x)$.
    \end{enumerate}
Moreover (note that, by definition, $\operatorname{dom}_{\text{\rm Fr}}(\nabla f)\subseteq \operatorname{dom}(\nabla f)$),
    \begin{enumerate}
        \item[(ii)] if $x\in \operatorname{dom}(\nabla f)$, then $y_n \rightharpoonup y=\nabla f(x)$;
        
        \item[(iii)] if $x\in \operatorname{dom}_{\text{\rm Fr}}(\nabla f)$, then $y_n \to y=\nabla f(x)$. 
    \end{enumerate}
    \end{Lemma}
\subsubsection{Hilbertian null sets}\label{sect:nullsets}
Before formally stating our Hilbertian version of McCann's theorem, we need  infinite-dimensional extensions of the finite-dimensional conditions of absolute continuity (i.e., ${\rm P}\in\Pac(\R^d)$) and Borel measures with Hausdorff dimension $(d-1)$ (i.e., ${\rm P}\in\PH(\R^d)$). Due to the absence of a Lebesgue measure, on general Hilbert spaces, this requires some care.

\begin{Definition}[Non-degenerate Gaussian distribution]
{\rm We say that a random variable $
\xi\in \H$ is  {\it non-degenerate Gaussian} if, for all $h \ne 0 \in \H$, the inner product $\langle
\xi , h\rangle\in \R$ is a non-degenerate Gaussian random variable, i.e., $\langle
\xi, h\rangle\sim \mathcal{N}(m_h, \sigma_h)$ with $\sigma_h>0$. The distribution $\mu_\xi$ of a non-degenerate Gaussian random variable $\xi$ is called a {\it non-degenerate Gaussian measure}. Denote by~$\operatorname{GN}(\H)$ the class of  Borel sets  negligible with respect to any non-degenerate Gaussian measure.
}
\end{Definition}

In the Euclidean space $\R^d$, the null sets of all nondegenerate Gaussian measures are exactly the same, and are equivalent to the Lebesgue negligible sets; for $\H= \R^d$, thus, $\operatorname{GN}(\H)$ reduces to the class of Lebesgue-null Borel sets.  This is no longer the case in infinite-dimensional $\H$, where several mutually singular non-degenerate Gaussian distributions exist (see e.g., the Feldman–H\'{a}jek theorem); this is why the definition of~$\operatorname{GN}(\H)$ imposes negligibility with respect to {\it any} non-degenerate Gaussian measure. 

Equivalently, $\operatorname{GN}(\H)$ can be described as the class of Borel sets that are negligible under any cube measure (see \cite{Csrnyei1999AronszajnNA}); a  cube measure is the distribution of  a random variable $a+\sum_{i=1}^n X_i e_i$ where  the span of $\{ e_i\}_{i\in \N}$ is dense in $\H$, such that $ \sum_{i\in\N} \|e_i\|^2<\infty$, and $\{X_i \}_{i\in \N}$ are uniformly distributed independent random variables with values in $(0,1)$. Moreover, \cite{Csrnyei1999AronszajnNA} proved that the  class of Aronszajn null sets previously mentioned  also coincides with {$\operatorname{GN}(\H)$}.

\begin{Definition}[Regular probability measures]\label{defn:Regular}
{\rm A probability measure ${\rm P}$ is called \emph{regular} if ${\rm P}(A)=0$ for all $A\in\operatorname{GN}(\H)$. In case $\H=\R^d$ for some finite $d$, the class of regular probability measures over~$\H$ coincides with  the class $\Pac(\R^d)$ of Lebesgue absolutely continuous measures, and we therefore denote  by $\Pac(\H)$ the family of all regular probability measures on $\H$.} 
\end{Definition}
Note that $\Pac(\H)$ contains all probability measures which are absolutely continuous with respect to some (degenerate or non-degenerate)  Gaussian measure, as well as all the Gaussian measures themselves. As we shall see, ${\rm P}\in\Pac(\H)$ is sufficient for existence and uniqueness in Theorem~\ref{Theorem:Maccan} below. But it is not necessary: as in the Euclidean case, where  it is sufficient for  ${\rm P}$ to be in the class~$\PH(\R^d)$ of measures giving mass zero to $(d-1)$-rectifiable\footnote{A set is called $(d-1)$-rectifiable if it can be written as a countable union of $\mathcal{C}^1$ manifolds, apart from a set of~$(d-1)$-dimensional Hausdorff measure zero  \citep[p.~271]{villani2008optimal}.} {sets}, this  assumption on ${\rm P}$  can be relaxed. 
Additionally, the  class $\PH({\mathbb R}^d)$  turns out to be too restrictive even in the Euclidean case---all we need is to ensure that the gradients of continuous l.s.c.~convex functions are~${\rm P}$-a.e.~well defined. 

Before  moving on with this discussion, let us formally introduce  the classes of probability measures we will need in this paper.

\begin{Definition}\label{defn:Lipschitz-Surfaces} 
{\rm (i) 
 A set~$ A_v\subseteq \H $, where $v\in \H\setminus\{0\}$, is called a {\em delta-convex  hypersurface}  
  if there exist two convex Lipschitz functions $ \tau_1,\tau_2: Z \to  \R$, with $Z \coloneqq \{ \lambda \, v\,: \ \lambda\in \R\}^{\perp} $ (the orthogonal complement of the space generated by $v$),  such that~$A_v = \{z+ (\tau_1(z)-\tau_2(z)) v:\  z\in Z \}$. Denote by $\Pdc(\H)$ the class of distributions giving mass zero to all delta-convex hypersurfaces.

 (ii) A   set of the form   $\{z+ \tau(z) v:\  z\in Z\}$, where  $\tau:\H\to \R$ is a  Lipschitz function, is called a  {\em Lipschitz
hypersurface}. Denote by $\Pl(\H)$ the class of distributions giving mass zero to all Lipschitz  hypersurfaces.
}
\end{Definition}
A delta-convex hypersurface is automatically  Lipschitz and   Lipschitz
hypersurfaces are Gaussian null sets:  hence, 
\begin{equation}  \label{contentionClasses}
\Pdc(\H) \supseteq \Pl(\H)\supseteq \Pac(\H)
\end{equation} 
(see e.g.,~\citet[p.~295]{Zajk1983} or~\citet[p.~521]{Zajk1978}).  The converse, however, is not true:~\citet[Example 1]{Zajk1979} shows that Lipschitz hyperspaces are not necessarily delta-convex hyperspaces, even in the Euclidean case.  
In the Euclidean case  ($\H =\R^d$) with $d\geq 2$, 
$$
\Pdc(\R ^d) \supsetneq \Pl(\R ^d)\supseteq \PH(\R ^d)\supsetneq \Pac(\R ^d),$$ where  $\PH(\R ^d) \subseteq \Pl(\R ^d)$ may be an equality. For $d=1$,   that is, for~$\H = \R$, however, we have~$\Pdc(\R) = \Pl(\R) = \PH(\R)=\Pac(\R )$. 

\begin{Remark}\label{remark:nonGauss}{\rm When $\H$ is infinite-dimensional, there exists a l.s.c.~convex function $f$ that is nowhere continuous whose gradient nevertheless pushes forward one non-degenerate Gaussian measure to another. Further, the set $\operatorname{dom}(\nabla f)$ is a Gaussian null set. We provide the construction of such a function below.
Let us consider a fixed orthonormal basis $\{e_i\}_{i\in \N}$ in the infinite-dimensional Hilbert space $\H$. Consider the unbounded operator $A:\operatorname{dom}(A)\to \H$ defined by $x\mapsto \sum_{i\in \N}4^i\langle x, e_i\rangle e_i \in \H$. Here, $\operatorname{dom}(A)$ is the pre-image~$A^{-1}(\H)$ of $\H$ under $A$,  that is, $\operatorname{dom}(A) =\{x\in\H: \ \sum_{i\in \N}8^i|\langle x, e_i\rangle|^2 \in~\!\R\}$. For the l.s.c.~convex function~$f: \H \to (-\infty, +\infty]$ defined as $f(x) \coloneqq \frac{1}{2} \|A x\|^2$ if $x \in \operatorname{dom}(A)$ and~$+\infty$ otherwise, the subdifferential is  $A x$ if~$x\in \operatorname{dom}(A)$, and is empty otherwise. Let $\{\xi_i\}_{i \in \N}\subset \R$ be a sequence of i.i.d.~$N(0,1)$  variables. We have the following two observations:}
\begin{enumerate}
    \item[(i)] The gradient of the l.s.c.~convex function $f$ pushes forward the Gaussian random variable~ $\sum_{i\in \N} \frac{1}{8^i} \xi_i e_i $ to $\sum_{i\in \N} \frac{1}{2^i} \xi_i e_i$. The function $f$ is discontinuous everywhere in $\H$. 

    \item[(ii)] The Gaussian random variable $X \coloneqq\sum_{i\in \N} \frac{1}{2^i} \xi_i e_i \sim~\mu$ is non-degenerate but  $$\mathbb{P}\left(X\in \operatorname{dom}(A)\right) = \mathbb{P}\left( \sum_{i\in \N} 2^i \xi_i^2 < +\infty \right)=0.$$  Since $\operatorname{dom}(\nabla f)=\operatorname{dom}(A)$,  however, $$\mu\left(\{ x\in \H: \ \partial f(x) \ \right.\left. \text{is a singleton}\}\right) = \mathbb{P}\left(X\in \operatorname{dom}(A)\right) = 0.$$ Thus, the set where the subdifferential is single-valued, i.e.,  $\operatorname{dom}(\nabla f)$, is a Gaussian null set. 
\end{enumerate}
\end{Remark}

\subsection{Existence and uniqueness of monotone measure-preserving
maps}\label{sec:Exist-Unique}
We can now state and prove our main result about the existence and uniqueness, without second-order moment restrictions, of monotone measure-preserving
maps in a separable Hilbert space $\H$---the Hilbertian generalization  of McCann's result in $\mathbb R^d$. 
\subsubsection{Existence and uniqueness}
\begin{Theorem}\label{Theorem:Maccan}
Let  ${\rm Q}\in \P(\H)$, where $\H$ is a separable Hilbert space. 
\begin{enumerate}
\item[(i)] If ${\rm P}\in 
\Pl(\H)$,  there exists a gradient of  convex function~$\nabla\psi$  pushing  ${\rm P} $ forward to ${\rm Q}$; 
\item[(ii)]if ${\rm P}\in 
 \Pl(\H)$ and $\operatorname{supp}({\rm Q})$ is bounded, then $\nabla\psi$ is unique ${\rm P}$-a.s.
\end{enumerate}
\end{Theorem}
The assumption of a boundedly supported $\rm Q$ is quite natural when the objective is the definition, without any moment assumptions,  of a Hilbertian transport-based notion of ``center-outward'' distribution or multivariate rank function similar to the finite-dimensional concepts proposed in \cite{chernoetal17} or \cite{Hallin2020DistributionAQ}:  the reference distributions  $\rm Q$ there, indeed, are (a) the Lebesgue uniform over the $d$-dimensional unit cube or (b) the spherical uniform over the unit ball in $\R ^d$. Natural infinite-dimensional extensions  would include (a) cubic probability measures, i.e., the distributions of  random variables of the form $\sum_{i\in \N} \lambda_i U_i e_i$ where $ \{\lambda_{i}\}_{i\in \N}\subset [0, +\infty)$ with~$\sum_{i\in \N} \lambda_i^2<+\infty$, $\{e_i\}_{i\in \N}\subset \H$ an orthonormal basis of $\H$, and $\{ U_i\}_{i\in \N}$ a sequence of i.i.d.\ univariate ${\rm Uniform}(0,1)$ variables (see e.g.,~\cite{Csrnyei1999AronszajnNA}) or (b) the distributions of  random variables of the form 
$U G/\|G\|$ where~$U\sim {\rm Uniform}(0,1)$ and $G$ is a Gaussian random variable in~$\H$. 

The proof of Theorem~\ref{Theorem:Maccan}, presented in Section~\ref{sec:McCann-proof}, relies on  four lemmas which we state and prove next.
\subsubsection{Four lemmas}
Throughout, it is tacitly assumed that $\H$ is a separable Hilbert space. 
\begin{Lemma}\label{Lemma:MCann}
Let $\{{\rm P}_n\}_{n\in \N}$ and $ \{{\rm Q}_n\}_{n\in \N}$ be sequences of probability measures in $\P(\H)$ such that~${\rm P}_n\xrightarrow{w} {\rm P}\in \P(\H)$ and ${\rm Q}_n\xrightarrow{w} {\rm Q}\in \P(\H) $. Suppose that the sequence $ \{ \gamma_n\}_{n\in \N}\in \P(\H\times \H) $ is such that $\gamma_n\in \Pi({\rm P}_n, {\rm Q}_n)$ and  
$ \operatorname{supp}(\gamma_n) $ is cyclically monotone for all $n\in \N$. Then, for any subsequence $\{\gamma_{n_k}\}_{k\in \N}$,  there exists a further subsequence $\{\gamma_{n_{k_i}}\}_{i\in \N}$ such that~$\gamma_{n_{k_i}}\xrightarrow{w}\gamma$ for some~$\gamma\in \Pi({\rm P}, {\rm Q})$ with cyclically monotone support. 
\end{Lemma}
\begin{proof}  Lemma 4.4 in  \citet{villani2008optimal} implies the tightness of  $\{ \gamma_n\}_{n\in\N}$. Hence, for any subsequence~$\{\gamma_{n_k}\}_{k\in \N}$,  there exists $\gamma$ and a further subsequence $\{\gamma_{n_{k_i}}\}_{i\in \N}$ such that $\gamma_{n_{k_i}}\xrightarrow{w}\gamma$. Let us  prove that $\gamma\in \Pi({\rm P}, {\rm Q})$ and that $\operatorname{supp}(\gamma) $ is cyclically monotone. For ease of notation, we keep the notation $\{ \gamma_n\}_{n}$ for the subsequence. 

Suppose that $\gamma$ is not cyclically monotone. Then, the  the subset 
\begin{align*}
M\coloneqq\Big\{ (x,y) &\text{ such that there exists a finite sequence }  \{(x_k,y_k)\}_k^n \text{ with } (x_1,y_1)=(x,y) \\
&\hspace{81mm}\text{ violating \eqref{Cyclically}   for }  y_{n+1}=y_1
\Big\}
\end{align*}
of $\H\times\H $  has strictly positive  $\gamma$ probability. 
That set $M$ is open, so that the Portmanteau theorem applies, yielding  
 $ \lim\inf_n\gamma_n\left(M\right)\geq \epsilon>0, $
which is impossible in view of the cyclical  monotonicity, for all $n\in \N$, of $\gamma_n$. Let $f:\H\to \R$ be a continuous bounded function. The function  $g(x,y)=f(x)$ is continuous and bounded  in $\H\times\H$, hence uniformly $\gamma_n$-integrable. Thus, $ \gamma_n(g)= {\rm P}_n(f)\to {\rm P}(f) $ 
and $\gamma_n(g)\to \gamma(g), $
so that the left  marginal of  $\gamma$ is  ${\rm P}$. Similarly,  ${\rm Q} $  is the right marginal.  
Any weak limit $\gamma$ of $\{\gamma_n\}_{n\in \N}$ has, thus,  cyclically monotone support and belongs to $\Pi({\rm P}, {\rm Q})$. The claim  follows. 
\end{proof}
Noting that the directional derivative and the subgradient $\partial f$ of a proper l.s.c.~convex\linebreak function~$f:\H \to (-\infty, \infty]$  are  related  through the formula  
\begin{equation}
    \label{limitsGateus}
    \lim_{t\to 0^+}\frac{f(h+t a)-f(h)}{t}=\sup_{y\in \partial f(h) }\langle y, a\rangle\qquad h,\, a\in \H
\end{equation}
(see, e.g.,~\citet[Theorem 17.19]{bauschkeMonotoneHilbert}), we obtain the following mean value theorem for convex functions in Hilbert spaces (the Hilbert space version of  \cite{McCann}'s  finite-dimensional Lemma 21). 
\begin{Lemma}\label{Lemma:MeanValue}
Let $f,g:\H\to (-\infty, +\infty]$ be  proper l.s.c.\   convex functions and $p,q\in{ \operatorname{cont}(f)\cap \operatorname{cont}(g)}$ be such that $ f(p)-f(q)=g(p)-g(q)$. Then, there exists $x_t=t p+(1-t)q$,\linebreak  with~$t\in (0,1)$, $u\in \partial f(x_t)$, and $v\in \partial g(x_t)$  such that 
$\langle u-v,p-q \rangle=0$.
\end{Lemma}
\begin{proof} 
By convexity, $f(x_t)$ and $g(x_t)$ are finite for any $x_t=t p+(1-t)q$ with $t\in (0,1)$. Moreover, the function $h:[0,1]\ni t\mapsto f(x_t)-g(x_t)\in \R$ is continuous (see~\citet[Corollary 9.20 and Theorem 8.29]{bauschkeMonotoneHilbert}). The values of $h(t)$  at $t=0$ and $t=1$ are the same by hypothesis; a \emph{fortiori} an {extreme value of $h$ in $[0,1]$} is attained at some $ t^*\in (0,1)$. Suppose, without loss of generality, that $h(t^*)$ is a maximum.   Letting $x^*\coloneqq t^* p+(1-t^*)q$, note that $f$ and $g$  both are continuous at $x^*$,   so that  $\partial f(x^*)$ and $\partial g(x^*)$   both are convex and weakly compact (see~\citet[Proposition 16.14]{bauschkeMonotoneHilbert}). Since the function $x\mapsto \langle x, p-q\rangle$ is weakly continuous, there exist  $u_+$, $u_-$, $v_+$, and $v_-$ such that 
$$   u_+=\arg\max_{y\in \partial f(x^*)}\langle y, p-q \rangle \quad \text{and} \quad u_-=\arg\min_{y\in \partial f(x^*)}\langle y, p-q \rangle, $$
for $f$, and 
\begin{align*}
    v_+=\arg\max_{y\in \partial g(x^*)}\langle y, p-q \rangle \quad \text{and} \quad v_-=\arg\min_{y\in \partial g(x^*)}\langle y, p-q \rangle
\end{align*}
for $g$. Thus, via \eqref{limitsGateus}, we obtain 
$$ \lim_{t\to 0^+}\frac{h(x^*+t{(p-q)})-h(x^*)}{t}= \langle u_+-v_+, p-q \rangle$$
and 
$${\lim_{t\to 0^-}\frac{h(x^*+t(p-q))-h(x^*)}{t}=  \lim_{t\to 0^+}\frac{h(x^*+t (q-p))-h(x^*)}{-t}=\langle u_--v_-, p-q \rangle.}$$
\citet[Lemma 20]{McCann} states that, as $h$ has a maximum at $t^*$, 
$$ \langle u_--v_-, p-q \rangle=\lim_{t\to 0^-}\frac{h(x^*+t a)-h(x^*)}{t}\leq 0\leq \lim_{t\to 0^+} \frac{h(x^*+t a)-h(x^*)}{t} =\langle u_+-v_+, p-q \rangle .$$
Hence, there exists some $\lambda \in [0,1]$ such that 
$\langle (u_--v_-)\lambda+(1-\lambda)(u_+-v_+), p-q \rangle=~\!0$.  Since~$\partial f(x^*)$ and $\partial g(x^*)$ are  compact, $ u_{\pm}$ and $v_{\pm}$ belong   to  $\partial f(x^*)$ and $\partial g(x^*)$, respectively. So, using the convexity of $\partial f(x^*)$ and $\partial g(x^*)$, we obtain $$ \lambda u_-+(1-\lambda)u_+\in \partial f\quad\text{and}\quad
 \lambda v_-+(1-\lambda)v_+\in~\!\partial g,$$ which concludes the proof.
\end{proof}

The next result states that if the gradients of two continuous  convex functions $f$ and $g$ differ at a point $p$, there exists a neighborhood of $p$ such that the set where both functions agree in this neighborhood is ``small'' (i.e., ${\rm P}$-negligible). Moreover, the inverse image by $\nabla g$ of $\partial f(\{f \ne g\})$ lies at a strictly positive distance from $p$.


\begin{Lemma}\label{Lemma:alexandrov}
Let $f$ and $g$ be  two proper l.s.c.\   convex functions from $\H$ to $(-\infty, +\infty]$ such   that, for some $p\in \operatorname{dom}(\nabla f)\cap \operatorname{dom}(\nabla g)\cap \operatorname{cont}(f) \cap \operatorname{cont}(g)$, $\nabla f(p)\neq \nabla g(p)=0$  and $ f(p)=g(p)=0$. Then, 
\begin{compactenum}
\item[(i)] $ {\cal X}\coloneqq (\nabla g)^{-1}(\partial f(\{h\in \H: \ f(h)>g(h) \}))\subseteq \{h\in \H: \ f(h)>g(h) \} $, 
\item[(ii)]$\inf_{h\in {\cal X}} \|h-p \|>0 $, and
\item[(iii)]there exists a neighborhood $\mathcal{U}_p \subset \H$ of $p$ such that the 
set~$\{h\in \mathcal{U}_p: \ f(h)=g(h) \} $  lies in a Lipschitz hypersurface.
\end{compactenum}
\end{Lemma}
\begin{proof}
The  proof is inspired by that of \citet[Lemma 13]{McCann}. Part $(i)$ of the lemma directly follows from the definition of subdifferentials. Take $ x\in {\cal X}$. Then, $ x\in  \operatorname{dom}(\nabla g)$ is such that~$\nabla g (x)\in \partial f(\{h\in \H: \ f(h)>g(h) \})$. Hence, 
\begin{equation}
    \label{rewritting}
    \nabla g (x)\in \partial f(h)\quad \text{  for some } h\in \operatorname{dom}(\partial f)\cap \{h\in \H: \ f(h)>g(h) \} 
\end{equation}
and thus, for any $z\in \H$,
\begin{equation}
\label{eq:LemmaAlex}
     f(z)-f(h)\geq \langle \nabla g (x) ,z-h\rangle \quad \text{and}\quad g(h)-g(x)\geq \langle \nabla g (x) ,h-x\rangle.
\end{equation}
In particular, for $x=z$, we obtain 
$$ f(x)-f(h)\geq \langle \nabla g (x) ,x-h\rangle \quad \text{ and }\quad g(h)-g(x)\geq \langle \nabla g (x) ,h-x\rangle.$$  
Since $ f(h)>g(h)$,  adding the above two inequalities yields 
$$ f(x) -g(x)> f(x)-f(h)+g(h)-g(x)\geq 0. $$
Hence $x\in \{h\in \H: \ f(h)>g(h) \}.$ This completes the proof of part $(i)$.

To prove part $(ii)$, let us assume that there exists a sequence $\{x_n\}_{n\in \N} \subseteq \cal X$ converging to $p$ in norm. Then $\nabla g(x_n)\in \partial f(\{h\in \H: \ f(h)>g(h) \}) $, so that there exists a sequence~$\{h_n\}_{n\in \N}\subseteq \{h\in \H: \ f(h)>g(h) \} $ such that $\nabla g(x_n)\in \partial f(h_n)$ for all $n\in \N$. Since $\nabla g(p)= 0$ and  $g(p)=0$, the  convexity of $g$ implies that $g\geq 0$. Moreover, the strong-to-weak continuity of the directional derivative (see Lemma~\ref{Lemma:strongtoweak}) implies $ \nabla g(x_n)\rightharpoonup \nabla g(p)=0$. It also follows from
\begin{equation}
    \label{eq:convOfG}
    |g(x_n)-g(x)|\leq (\|\nabla g(x_n)\|+\|\nabla g(x)\|)\|x_n-p\|
\end{equation}
that~$ g(x_n)\to   g(p)=0$. On the other hand, $\nabla f(p)\neq 0$ and the fact that, for all $z\in \operatorname{dom}(\nabla f)$, 
$$-f(z)= f(p)-f(z)\geq \langle \nabla f (z) ,p-z\rangle $$
jointly imply (taking $z_n \coloneqq p-\frac{1}{n}\nabla f (p)$, with $n\in \N$ large enough) that  
$$f(z_n)\leq -\frac{1}{n}\langle \nabla f (z_n) ,\nabla f (p)\rangle. $$
Hence, due to the strong-to-weak continuity of the directional derivative, 
$${\lim\sup_{n \to \infty} n\, f(z_n)\leq\! -\|\nabla f (p) \|^2  },$$ and  there exists $z=z_{n_0}$ (for $n_0$ big enough) such that $f(z)<0$. Using \eqref{eq:LemmaAlex}, we obtain 
\begin{equation*}
     (f(z)-f(h))+(g(h)-g(x))\geq \langle \nabla g (x) ,z-h+(h-x)\rangle=\langle \nabla g (x) ,z-x\rangle,
\end{equation*}
for any $z\in \mathcal{H}$ and $h,x$ as in \eqref{rewritting}. Since $ f(h)>g(h)$, we have 
$
     f(z)-g(x)\geq \langle \nabla g (x) ,z-x\rangle.
$  By taking $z = z_{n_0}$ and $x = x_n$, 
this yields  
\begin{align}
    0>f(z)&\geq \langle \nabla g(x_n), z-x_n \rangle +g(x_n) \nonumber
    \\ 
    &\geq \langle \nabla g(x_n), z \rangle-\langle \nabla g(x_n), x_n-p\rangle-\langle \nabla g(x_n), p\rangle +g(x_n) \nonumber 
    \\
    &\geq \langle \nabla g(x_n), z-p\rangle-\| \nabla g(x_n)\|\| x_n-p \| +g(x_n).
\label{last}\end{align}
The right-hand side in \eqref{last} tends to $0$ since: (i) the first term goes to zero as $ \nabla g(x_n)\rightharpoonup \nabla g(p)=~\! 0$; (ii) the second term goes to zero by the boundedness of $\| \nabla g(x_n)\|$ (see, e.g.,~\citet[Proposition~3.13 (iii)]{Brezis}) and the fact that $x_n\to p$; and (iii) the last term tends to 0 by \eqref{eq:convOfG}.
This, however, yields a contradiction from which  we conclude that  $p\not\in \overline{\cal X}$. This completes the proof of part $(ii)$ of the lemma.

Turning to part $(iii)$, let us write $\{f=g\}$ for $\{x\in \H: \ f(x)=g(x) \} $. {Consider an orthonormal basis $\{ e_i\}_{i\in \N}$ such that~$e_1\coloneqq  {\|\nabla f(p) \| ^{-1}}\nabla f(p)\neq 0$.}  
As in \cite{McCann}, the goal is to show the existence of a Lipschitz function $\tau:\H \to \R$ and a neighborhood $\mathcal{U}_p$ of $p$  such that 
\begin{equation}\label{Lipchitz:rectificable}
    \{f=g\}\cap \mathcal{U}_p\subseteq \{z + \tau(z) e_1: z \in Z\} 
\end{equation}
where $Z = \{\lambda e_1: \lambda \in \R\}^\perp\vspace{0mm}$ and $\pi(x)\coloneqq x- \langle e_1, x\rangle e_1$ is the orthogonal projection of $x \in \H$ onto   the closure 
of the subspace generated by $\{e_i\}_{i=2}^{\infty}$. 
{ Set $h\in\operatorname{cont}(f)\cap \operatorname{dom}(\nabla f)$.}
The strong-to-weak continuity of the subdifferential (see, e.g.,~\citet[Proposition 17.3]{bauschkeMonotoneHilbert}) means that $u_n\rightharpoonup \nabla f(p)= :\lambda_1 e_1$ and~$v_n\rightharpoonup \nabla g(p)$ for any $p_n\to p$,  $u_n\in \partial f(p_n)$ and $v_n\in \partial g(p_n)$.  Let~$\mathcal{U}_p$ denote a neighborhood of $p$ small enough  that {$\partial f(\mathcal{U}_p)\cup \partial g(\mathcal{U}_p)$} is bounded,  i.e.,  contained in a ball $\mathcal{B}(0,R/2)$ (see, e.g., \citet[Pro\-po\-si\-tion~16.14]{bauschkeMonotoneHilbert}). Then we can find a neighborhood (we use the same notation $\mathcal{U}_p$ as above) of $p$ such that 
\begin{equation}\label{eq:Lemma:Lip}
    \langle e_1, u-v\rangle >\frac{\lambda_1}{2} \qquad\text{and} \qquad  \| \pi(u)-\pi(v)\|\leq R
\end{equation}
whenever  $u\in \partial f(x)$, $v\in \partial g(x)$,  and~$x\in \mathcal{U}_p$. The first inequality in \eqref{eq:Lemma:Lip} holds due to the strong-to-weak continuity of the subdifferential (namely, for $x \approx p$, we have $ \langle e_1, u \rangle \approx \langle e_1, \nabla f(p) \rangle = \lambda_1$, and~$\langle e_1, v \rangle \approx \langle e_1, \nabla g(p) \rangle = 0$). As for the second inequality in \eqref{eq:Lemma:Lip}, noting that  the projection operator $\pi$ is 1-Lipschitz, it  follows from the fact that  $u,v \in \mathcal{B}(0,R/2)$.

We now proceed with the construction of the Lipschitz function $\tau$. It follows from~\eqref{eq:Lemma:Lip}  that~$ \langle e_1, u-v\rangle>0$. For $x \in \H$ and $t \in \R$ such that $p + \pi(x)+t\, e_1\in \mathcal{U}_p$,  define the real functions $h_{x}(t)\coloneqq[f-g](p+\pi(x)+t\, e_1)$.  Observe that $t\mapsto h_{x}(t)$ is strictly monotone   in its domain. To see this, suppose that for $t_1 \ne t_2$, we have $h_{x}(t_1) = h_{x}(t_2)$; then by Lemma~\ref{Lemma:MeanValue} we would have $(t_1-t_2) \langle e_1, u -  v \rangle = 0$ which, letting~$s = t^* t_1 + (1- t^*) t_2$ with~$t^* \in (0,1)$, is a contradiction by~\eqref{eq:Lemma:Lip} (here $u \in \partial f(p+\pi(x)+ s e_1)$ and $u \in \partial g(p+\pi(x)+s e_1)$.

Let $t \ne 0$ be such that  $p+t \, e_1, p-t \, e_1\in \mathcal{U}_p$. We can pick $r>0$ such that $$h_0(-t)<-2r<0 = h_0(0) <2r<h_0(t)$$  and, by the continuity of $f-g $, also~$h_x(-s)<-r<0<r<h_x(s)$ for all $p+\pi(x)-s\, e_1$\linebreak and~$p+\pi(x)+s\, e_1$ belonging, respectively,  to  balls, $\mathcal{B}(p-t \, e_1, \rho)$ and $\mathcal{B}(p+t\, e_1, \rho)$, say, included in a small neighborhood of $p$ ensuring that $s$ and $\|\pi(x)\|$ are small. We can assume that such neighbor- hoods~are contained in $\mathcal{U}_p$. The intersection of $\mathcal{U}_p$ and the cylinder $\{x\in \H:\ \| \pi(x-p)\|<~\!\rho\}  $  is still a neighborhood of $p$ that we still denote as  $\mathcal{U}_p$. 

Let $\mathcal{V}\coloneqq \{x\in \H:p + \pi(x) \in \mathcal{U}_p\}$. For any $x \in \mathcal{V}$, $t\mapsto h_{x}(t)$ is a strictly monotone and continuous function taking both positive and negative values  in a neighborhood of 0. Then there exists a unique~$t_x$ such that $h_{x}(t_x)=0$. By construction, $t_x$ depends only on $\pi(x)$: writing $t_{\pi(x)}$ instead of $t_x$, define~$\tau: \mathcal{V} \to \R$ such that $\tau:x \mapsto t_{\pi(x)}$. Thus, $\tau(x) = \tau(\pi(x))$ if $\pi(x) \in \mathcal{V}$.
 Let us show that $\tau$ is Lipschitz. Set~$x, w\in \mathcal{V}$. Since $$ 0=h_x(t_{\pi(x)}) = [f-g](p+\pi(x)+t_{\pi(x)}\, e_1)=[f-g](p+\pi(w)+t_{\pi(w)}\, e_1) = h_w(t_{\pi(w)}),$$ Lemma~\ref{Lemma:MeanValue} ensures the existence of some $u\in \partial f(y)$ and $v\in \partial g(y)$ where $$y \coloneqq p + \lambda \pi(x) + (1-\lambda) \pi(w) + \lambda t_{\pi(x)} + (1-\lambda) t_{\pi(w)}$$ for some $\lambda \in (0,1)$ such that $\langle u-v, \pi(w-x)+(t_{\pi(x)}-t_{\pi(w)})\, e_1\rangle=0 $. Observe that $y \in \mathcal{U}_p$.  This further implies that 
 $$\langle u-v, \pi(w-x) \rangle = - (t_{\pi(x)}-t_{\pi(w)})\, \langle u-v, e_1\rangle , \text{ hence }|\langle u-v, \pi(w-x) \rangle| = |\tau(x)- \tau(w)| |\langle u-v, e_1\rangle|.$$  
In view of \eqref{eq:Lemma:Lip}, we thus obtain 
\begin{eqnarray*}
\frac{ \lambda_1}{2} | \tau(x)-\tau(w)|  &<& | \tau(x)-\tau(w)| \,| \langle u-v, e_1\rangle|
= |\langle u-v, \pi(w-x)\rangle| \\
& =& |\langle \pi(u-v), \pi(w-x)\rangle|  
 \le  \| \pi(u-v) \|\, \|\pi(w-x)\| \\
& \le &  R\, \|\pi(w-x)\|.
\end{eqnarray*}
Then $\tau$ is $(2R/\lambda_1)$-Lipschitz. Note that such a $\tau$ can be extended to the whole space $\H$ while preserving the Lipschitz constant and the dependence  of $\tau(x)$   on $\pi(x)$ only. To prove this, we only need to apply  \cite[Theorem~1]{HIRIARTURRUTY1980539} to the restriction of $\tau(\cdot)$  to the set $\pi(\H)$.  For any $x \in \H$, let us define the {\it translation of} $\tau$ as $\tau^*(p + x) \coloneqq \tau(x) = \tau(\pi(x))$ and show that~$\tau^*$ satisfies \eqref{Lipchitz:rectificable}.  Take $h \in \{f = g\} \cap \mathcal{U}_p$. Then, $ h = p + \pi(x) + s e_1$, for some $s \in \R$ and~$x \coloneqq h-p$. Note that $\tau(x)$ is the unique point in the line $p + \pi(x) + t e_1$ such that, for $t$ small,~$f(p + \pi(x) + \tau(x) e_1) = g(p + \pi(x) + \tau(x) e_1)$. But $h \in \{f = g\} \cap \mathcal{U}_p$ is also a point in the line segment 
$p + \pi(x) + t e_1$ such that, for $t$ small,   $f(h) = g(h)$. By the uniqueness of the construction of $\tau$, we get $\tau(x) = s$. Therefore, $h = p + \pi(x) + \tau(x) e_1 = z + \tau^*(z) e_1$ where~$z \coloneqq p + \pi(x)$ (note that $\tau^*(z) = \tau(\pi(x)) = \tau(x)$). Thus,~\eqref{Lipchitz:rectificable} holds for $\tau = \tau^*$, which completes the proof.
\end{proof}

\subsubsection{Proof of Theorem~\ref{Theorem:Maccan}}\label{sec:McCann-proof}

We now turn to the proof of Theorem~\ref{Theorem:Maccan}. We first assume that ${\rm Q}$ has  bounded support and, under this assumption, we prove $(i)$   existence  and $(ii)$ uniqueness. In $(iii)$, we then extend the existence result to the case when $\operatorname{supp}({\rm Q})$ is unbounded.

\noindent{$(i)$ [Existence, boundedly supported ${\rm Q}$.]} 
This part of the proof follows along similar steps as in  the finite-dimensional case (cf.~\citet[Theorem 6]{McCann})---with the significant difference, however, that the Riesz-Markov  theorem no longer can be invoked since Hilbert spaces are not necessarily  locally compact: the space of Radon measures and   the dual of the space of bounded continuous functions, thus, do not necessarily coincide. 



We can easily construct two sequences of probability measures $ {\rm P}_n$ and ${\rm Q}_n$  with finite second-order moments for all $n$ converging weakly to~${\rm P}$ and ${\rm Q}$, respectively. The existence of measure-preserving cyclically monotone maps between~$ {\rm P}_n$ and ${\rm Q}_n$ follows from \cite{Cuesta1989NotesOT}. Therefore, we can construct a sequence $\gamma_n\in \Pi({\rm P}_n, {\rm Q}_n)$ with 
$ \operatorname{supp}(\gamma_n) $ cyclically monotone.  By  Lemma~\ref{Lemma:MCann} and \cite{RockafellarMaximalMonot}, there exists  $\gamma\in \Pi({\rm P}, {\rm Q})$ and a proper l.s.c.\ convex function $\psi$ such that $\partial\psi\supseteq\operatorname{supp}(\gamma)$. Note that, denoting by $\psi^*$   the convex conjugate of  $\psi$ and by~$\bar{\operatorname{ch}}(\operatorname{supp}({\rm Q}))$ the closed convex hull of $\operatorname{supp}({\rm Q})$, 
\begin{equation}\label{eq:countinous}
\bar{\psi}(x)\coloneqq\sup_{y\in\bar{\operatorname{ch}}(\operatorname{supp}({\rm Q}))}\{ \langle x,y \rangle -\psi^*(y)\}
\end{equation}
and $\psi$ {\rm P}-a.e.~agree (see~\citet[p.~147]{AmbrosioGradient}). For $x,x' \in \H$, as 
$$ |\bar{\psi}(x)-\bar{\psi}(x')|\leq \sup_{y\in\bar{\operatorname{ch}}(\operatorname{supp}({\rm Q}))}\{ \langle x-x',y \rangle\}\leq \|x-x'\|\sup_{y\in\bar{\operatorname{ch}}(\operatorname{supp}({\rm Q}))}\|y\|,
$$
the function $\bar{\psi}$ is continuous. Accordingly, we keep the notation $\psi=\bar{\psi}$.

Since~$\rm P \in \Pl(\H)$,  $\partial\psi(x)$ is a singleton for ${\rm P}$-a.e.~$x$, i.e.,  ${\rm P}(\operatorname{dom}(\nabla\psi)) = 1$ (see~\cite{Zajk1979}). Define $$T: \operatorname{dom}(\nabla\psi) \to \H \;\; \mbox{with} \;\;  x\mapsto T(x)= \nabla \psi(x)\in  \H. $$
This $T$ is a Borel map and is such that $\gamma = ({\rm identity} \times T){\#}{\rm P}$ (as $\partial\psi\supseteq\operatorname{supp}(\gamma)$). 
Thus $T = \nabla \psi$ is the gradient of a convex function  pushing ${\rm P}$ to ${\rm Q}$, thereby completing the proof of the  existence part of Theorem \ref{Theorem:Maccan} for boundedly supported $\rm Q$. \\

\noindent {$(ii)$ [Uniqueness, boundedly supported ${\rm Q}$.]} To prove  uniqueness, let us assume that two  distinct l.s.c.\ proper convex functions~$f$ and $g$ are such that $\nabla f$ and $\nabla g$ both push ${\rm P} $ forward to ${\rm Q}$. 
To start with, assume  that~$p\in \operatorname{supp}({\rm P})\cap \operatorname{dom}(\nabla f)\cap \operatorname{dom}(\nabla g)\cap \operatorname{cont}(f)\cap \operatorname{cont}(g)$ is such that $\nabla f(p)\neq \nabla g(p) $. Consider the  two functions (from $\H$ to $\R$)
\begin{equation}\label{eq:f*}
f^*(\cdot )\coloneqq f(\cdot )-\langle\cdot, \nabla g(p) \rangle -(f(p)-\langle p, \nabla g(p)\rangle)
\end{equation}
and  
\begin{equation}\label{eq:g*}
g^*(\cdot )\coloneqq g(\cdot )-\langle\cdot, \nabla g(p)\rangle-(g(p)-\langle p, \nabla g(p)\rangle)
\end{equation}
obtained by adding to $f$ and $g$ affine functions (depending on $f(p)$, $g(p)$, and $\nabla g(p)$). Clearly,  one has~
$ f^*(p)=g^*(p)=0$ and ~$f^*$ and $g^*$  satisfy $\nabla f^*(p)\neq \nabla g^*(p)=0$ just as $f$ and $g$. For simplicity,  we keep  for~$f^*$ and $g^*$ the notation $f$ and $g$. 

Let $\mathcal{U}_p$ be the open neighborhood  of~$p$ given by Lemma~\ref{Lemma:alexandrov}.\footnote{To meet the  assumptions of Lemma~\ref{Lemma:alexandrov}, we need to ensure that,    any  proper l.s.c.~convex function $f$ such that~$\nabla f$  pushes~${\rm P} $ forward to ${\rm Q}$ satisfies ${\rm P}(\operatorname{cont}(f))=1$. This, however, follows from~\citet[Proposition~17.41]{bauschkeMonotoneHilbert} and \eqref{eq:countinous}.} 
Then the set $$\mathcal{A} \coloneqq \mathcal{U}_p\cap\{h\in \H: \ f(h)= g(h) \}$$ is contained in a Lipschitz hypersurface so that, since $p\in \operatorname{supp}({\rm P})$ and ${\rm P}\in 
 \Pl(\H)\supsetneq \Pac(\H)$, we obtain  ${\rm P}(\mathcal{U}_p\setminus \mathcal{A})>0 $. Now,   $\mathcal{U}_p\setminus \mathcal{A}$ is the disjoint union of
 $$ M\coloneqq \mathcal{U}_p\cap\{h\in \H: \ f(h)> g(h) \}\quad \text{ and } \quad N\coloneqq \mathcal{U}_p\cap\{h\in \H: \ f(h)<g(h) \},$$  one  of which at least has strictly positive ${\rm P}$-probability; let us assume that ${\rm P}(M)>0$. Lemma~\ref{Lemma:alexandrov} implies that~$ {\cal X} \coloneqq(\nabla g)^{-1}(\partial f( M))\subseteq M $ and $\inf_{h\in \cal X} \|h-p \|>0 $, so that there exists an open neighborhood~$\mathcal{W}_p$  of $p$ such that $ {\cal X}\cap\mathcal{W}_p=\emptyset $. Assume, without loss of generality, that $\mathcal{W}_p\subseteq \mathcal{U}_p$.  As a consequence,~${\rm P}({\cal X}\cap (\H\setminus\mathcal{W}_p))\leq{\rm P}(M\cap (\H\setminus\mathcal{W}_p))$ and, by the previous argument, ${\rm P}(M\cap \mathcal{W}_p)>0$. Therefore, since $ {\cal X}\cap\mathcal{W}_p=\emptyset $,
\begin{equation*}
    {\rm P}({\cal X})={\rm P}({\cal X}\cap \mathcal{W}_p)+{\rm P}({\cal X}\cap (\H\setminus\mathcal{W}_p))<{\rm P}(M\cap \mathcal{W}_p)+{\rm P}(M\cap (\H\setminus\mathcal{W}_p))={\rm P}(M),
\end{equation*}
while, on the other hand, as $\operatorname{dom}(\nabla f)$ has ${\rm P}$-probability 1 (by part $(i)$ above),
$$ {\rm P}((\nabla f)^{-1}(\partial f( M)))= {\rm P}((\nabla f)^{-1}(\partial f( M))\cap \operatorname{dom}(\nabla f))\geq {\rm P}(M),$$
so that ${\rm P}((\nabla f)^{-1}(\partial f( M)))\neq {\rm P}({\cal X}) 
 = {\rm P}((\nabla g)^{-1}(\partial f( M)))$. This   contradicts the fact that both~$\nabla f$ and~$\nabla g$ (up to a translation: see~\eqref{eq:f*} and~\eqref{eq:g*}) are pushing   ${\rm P}$ forward to ${\rm Q}$. As a consequence, $\nabla f$ and $\nabla g$ must agree on $ \operatorname{supp}({\rm P})\cap \operatorname{dom}(\nabla f)\cap \operatorname{dom}(\nabla g)$. Uniqueness follows.\\  

\noindent $(iii)$ [Unbounded $\operatorname{supp}({\rm Q})$.]
The rest of the proof is inspired by that of \citet[Theorem 6.2.10]{AmbrosioGradient}. We have shown before that there exists  $\gamma\in \Pi({\rm P}, {\rm Q})$ and a proper l.s.c.\ convex function $\psi$ such that $\partial\psi\supseteq\operatorname{supp}(\gamma)$. Denoting by $\mathbb{I}_{\mathcal{B}(0,n)}$ the indicator function of the ball $\mathcal{B}(0,n)$ with radius $n$ centered at $0 \in \H$,   characterize the measure $\gamma_n$ as satisfying 
$$\int f(x,y)d\gamma_n(x,y)=\int \mathbb{I}_{\mathcal{B}(0,n)}(y)f(x,y)d\gamma(x,y)$$ for any continuous bounded function $f:\H\times \H\to \R$.  Since $\operatorname{supp}(\gamma)$ is cyclically monotone, $\operatorname{supp}(\gamma_n)$ is  cyclically monotone as well. Denote by ${\rm P}_n$ and ${\rm Q}_n$ the marginals of $\gamma_n$. Note that~$\int f(y)d{\rm Q}_n(y)=\int \mathbb{I}_{\mathcal{B}(0,n)}(y) f(y)d{\rm Q}(y)$ for any continuous bounded function $f:\H\to \R$, that~${\rm Q}_n$ is boundedly supported, and that~${\rm P}_n$ is absolutely continuous with respect to $\rm P$.  It follows from part $(i)$ of this proof (applied to the duly rescaled $\gamma_n$ to make it a probability measure) that~$\gamma_n$ is unique and that there exists a unique gradient $\nabla\psi_n$ of a l.s.c.~convex function~$\psi_n$ such\linebreak that~$({\rm Identity} \times \nabla \psi_n) {\#} {\rm P}_n=\gamma_n$. Since $\operatorname{supp}(\gamma_n)\subset \operatorname{supp}(\gamma)$ (and we know that $\partial\psi\supseteq\operatorname{supp}(\gamma)$), it follows from   parts $(i)$ and $(ii)$ of the proof above that $ \nabla \psi_n= \nabla \psi $, ${\rm P}_n$-a.s. 
As a consequence,~$\gamma_n=({\rm Identity} \times \nabla \psi) {\#} {\rm P}_n$, so that by taking weak limits as $n\to +\infty$, we\linebreak obtain~$\gamma=({\rm Identity} \times \nabla \psi) {\#} {\rm P}$. Such a $\nabla \psi$ satisfies the assumptions of Theorem~\ref{Theorem:Maccan}~(i) and, moreover, any $\gamma\in \Pi({\rm P}, {\rm Q})$ such that $\partial\psi\supseteq\operatorname{supp}(\gamma)$ for some proper l.s.c.\ convex function $\psi$ is of the form $\gamma=({\rm Identity} \times \nabla \psi) {\#} {\rm P}$. \qed

\section{Stability of  transport maps and a central limit theorem}\label{sec:stability}
The objective of this section is to establish the stability of transport maps of the form $\nabla\psi$,\linebreak where~$\psi :\H \to (-\infty,\infty]$ is a proper l.s.c.~convex function. 
This problem has not been addressed so far in the literature for infinite-dimensional spaces. Indeed, the techniques of \cite{Hallin2020DistributionAQ},~\cite{GhosalSenAOS},~\cite{Barrio2022NonparametricMC} or \cite{Segers2022GraphicalAU} are based on the Fell topology, which does not have nice properties in non-locally compact spaces. Consequently, as stressed by \citet[p.~4]{Segers2022GraphicalAU}, stability results, in general Hilbert spaces, are a challenging topic. 

Our main stability result is stated in Theorem~\ref{Theorem stability}, along with some remarks. In Section~\ref{CLTsec}, we establish a central limit result for the cost of the optimal transport. The proof of Theorem~\ref{Theorem stability}  is given in Section~\ref{Theorem stabilitySec}. 

\subsection{A stability result for cyclically monotone transport maps}\label{sec31}

Throughout, $\cal H$ stands for a separable Hilbert space. Unless otherwise specified,  limits are taken as~$n\to\infty$. By a {\it strongly compact set $K \subset \H$} we mean a compact set $K$ with respect to the strong topology (i.e., $x_n \to x$ if and only if $\|x_n - x\| \to 0$). Recall that, in a Hilbert space~$\H$, the closed unit ball $\mathcal{B}(0,1) = \{x \in \H: \|x\| \le 1\}$ is weakly compact, but may not be strongly compact; see, e.g.,~\citet[Theorem 3.16]{Brezis}.

\begin{Theorem}\label{Theorem stability}
    Let $\{{\rm P}_n\}_{n\in \N}$ and $\{{\rm Q}_n\}_{n\in \N}\subseteq\P (\H)$  be two sequences of probability measures such that ${\rm P}_n\xrightarrow{w}{\rm P}\in \Pl(\H)$ and ${\rm Q}_n\xrightarrow{w}{\rm Q}\in \P(\H)$. Assume that $\operatorname{supp}({\rm Q})$ and $ \operatorname{supp}({\rm Q}_n)$ both are  subsets of the ball $\mathcal{B}(0,M)$ for all $n\in \N$ and some $M>0$.  Let $ \gamma_n\in \Pi({\rm P}_n,{\rm Q}_n)$  be 
    such that~$\operatorname{supp}(\gamma_n)\subseteq  \partial \psi_n  $ for some proper l.s.c.~convex function $\psi_n$, and let $ \psi $ be a proper l.s.c.\ convex function such that $ \nabla \psi\# {\rm P}={\rm Q}$.
    Then, 
    \begin{itemize}
        \item[(i)] for any strongly compact set $K\subseteq\operatorname{int}(  \operatorname{dom}(\nabla \psi))\cap \operatorname{int}(\operatorname{supp}({\rm P}))$ and any $h\in \H$,
\begin{equation}\label{(i)}
     \sup_{(x,y)\in \partial \psi_n, x\in K } \langle y-\nabla\psi(x), h \rangle  \longrightarrow 0;
     \end{equation}
      \item[(ii)]  if $ \psi$ is (up to additive constants) the unique proper l.s.c.\ convex function such that $ \nabla \psi\# {\rm P}={\rm Q}$,  there exists a sequence $ \{ a_n\}_{n\in \N}$ such that, for any strongly compact convex set $K \subseteq \operatorname{supp}({\rm P})$, 
   \begin{equation}\label{(ii)} \sup_{x\in K}\left|\psi_n(x)+a_n-\psi(x)\right|\longrightarrow 0;\end{equation}
   \item[(iii)] if  $\operatorname{supp}(Q)$ is strongly compact or $\operatorname{dim}(\H)$ is finite, under the same assumptions as in $(i)$, we have
     \begin{equation}\label{(iii)} \sup_{(x,y)\in \partial \psi_n, x\in K } \left\| y-\nabla\psi(x)\right\|  \longrightarrow 0.\end{equation} 
    \end{itemize}

\end{Theorem}



{\rm Relaxing some of the assumptions in Theorem~\ref{Theorem stability} would be quite desirable; whether this is possible, however, is unclear. Below are two examples of violations of these assumption under which the theorem no longer holds true. } 
\begin{itemize}
\item[(a)]  {\rm [Boundedness of ${\rm supp}({\rm Q}_n)$ and ${\rm supp}({\rm Q})$] Let  ${\rm P}={\rm Q}$ denote the centered Gaussian distribution with covariance operator ${\rm diag}(1/4^i)$---namely, the distribution of $X = \sum_{i \ge 1} (1/2^i) \xi_i e_i$ where $\{\xi_i\}_{i\ge 1}$ are i.i.d.~$N(0,1)$ and $\{e_i\}_{i \in \mathbb{N}}$ is some fixed orthonormal basis of $\H$. \linebreak  Let ${\rm Q}_n$,   $n\in \N$, denote the centered Gaussian distribution with  covariance opera-\linebreak tor~${\rm diag}(1/(2-1/n)^{2 i })$ with respect to the same orthonormal basis $\{ e_{i}\}_{i\in \N}$. Assume that~$ {\rm P}_n={\rm P}$ for all $n\in \N$.  Clearly:   (i) the identity map ${\rm Id}$ is the unique gradient of a convex function that pushes ${\rm P}$ forward to ${\rm Q}$; (ii) 
the map $x \mapsto T_n (x)$ defined as
$$T_n(\sum_{i \ge 1}\lambda_i e_i) \coloneqq \sum_{i \ge 1}(2/(2-1/n))^{ i }\lambda_i e_i$$ 
is the unique gradient of a convex function that pushes ${\rm P}_n$ forward to ${\rm Q}_n$ (see e.g.,~\citet[Proposition~2.2]{CuestaAlbertos1996OnLB}); and (iii) we have $ {\rm P}_n\xrightarrow{w}{\rm P}$ and $ {\rm Q}_n\xrightarrow{w}{\rm Q}$. Let~$\{x_m\}_{m \ge 1}$ be a sequence such that $x_m \to x\coloneqq \sum_{i \ge 1} (1/i) e_{i} $; then $ \|T_n(x_m)\|\to \infty$ as $m\to \infty$ (where~$n$ is kept fixed). Thus, we can choose $x_{m(n)}$ such that $\|T_n(x_{m(n)})\|>n$. Set $h\in \H$. As the sequence $\{x_{m(n)}\}_{n\in \N}$ converges strongly to $x$, we have $\langle h, {\rm Id}(x_{m(n)})\rangle \to \langle h, x\rangle$ as $n \to \infty$. However,   $\|T_n(x_{m(n)})\|\to \infty$. This means that the sequence $\{T_n(x_{m(n)})\}_{n\in \N}$ is unbounded. Banach–Steinhaus's theorem (e.g., \citet[Theorem~2.2]{Brezis}) then implies the existence of some~$h\in \H$ such that $\langle h, T_n(x_{m(n)})\rangle\to \infty$ as $n \to \infty$. As a consequence,   $$|\langle h, T_n(x_{m(n)})-{\rm Id}(x_{m(n)})\rangle|= |\langle h, T_n(x_{m(n)})-x_{m(n)}\rangle |\to +\infty$$
as~$n \to \infty$, where we have used the fact that $\langle h, x_{m(n)}\rangle \to \langle h, x \rangle$. Hence, part $(i)$ of Theorem~\ref{Theorem stability} no longer holds true. 
}

\item[(b)]  {\rm   [$K\subseteq \operatorname{int}(\operatorname{supp}({\rm P}))$] This assumption also appears in the finite-dimensional case (see \cite{delbarrio2021central,Segers2022GraphicalAU,gonzalezdelgado2021twosample}). The following example, where $\H = \R$,  shows that this assumption  cannot be avoided. Let ${\rm P}={\rm Q}\in \Pac(\R)$ be the uniform distribution on $(0,1) \cup (2,3)$. Here, the identity function ${\rm Id}$ is the monotone transport map pushing ${\rm P}$ forward to ${\rm Q}$. Of course, ${\rm Id}$ is  everywhere single-valued. Let ${\rm P}_n$ be the uniform distribution on $(0,1+1/n) \cup (2+1/n,3)$. Then, a subdifferential $\partial \psi_n$ (defined as in Theorem~\ref{Theorem stability}) pushing ${\rm P}_n$ forward to ${\rm Q}$ is
\begin{equation*}
  \partial \psi_n(x) \coloneqq\left\lbrace  \begin{array}{ll}
         \{ x\}  &  {\rm if}\  x\in (-\infty ,1)\cup (2+1/n,+\infty) ,  \\
         \{ x+1\}  &  {\rm if}\  x\in (1+1/n),\\
        \{2+1/n\}  &  {\rm if}\  x\in [1+1/n, 2], \\
         \left[1,2\right]&  {\rm if}\  x=1 .
    \end{array} \right.
\end{equation*}
Clearly, 
$2\!\in\!   \partial \psi_n(1)$ for all $n\in \N$, but $2\not\to 1={\rm Id}(1)$.  This counterexample   arises from the fact that the identity function ${\rm Id}$ is not the unique monotone mapping  pushing  ${\rm P}$ forward to~${\rm Q}$. Although any other mapping $T$   satisfying this property agrees with ${\rm Id}$ in the interior of the support of ${\rm P}$, there exist mappings that differ from ${\rm Id}$ on the boundary of the supp(${\rm P}$). For example, consider the mapping $T(x)={2}$ for $x\in [1,2)$ and $T(x)=x$ otherwise. This mapping is monotone and does not agree with ${\rm Id}$ at   point $1\!\in\! \operatorname{supp}({\rm P})$. Hence, we conclude that the hypothesis $K\subseteq \operatorname{int}(\operatorname{supp}({\rm P}))$ cannot be relaxed without imposing some additional assumptions (such as strict convexity; see \cite{Segers2022GraphicalAU}) on the shape of the support of ${\rm Q}$. 
}
\item[(c)] [Strong compactness of $\operatorname{supp}(\rm Q)$] In general, the subdifferential of a l.s.c.~convex function is strong-to-weak continuous  in the set of differentiability points and strong-to-strong continuous in the set of Fr\'{e}chet differentiability points~\cite[Theorem~21.22]{bauschkeMonotoneHilbert}. The latter is not necessarily a null set with respect to any non-degenerate Gaussian measure (see \citet[Problem 5.12.23]{Bogachev2008GaussianM} for an example). Hence, in an infinite-dimensional space, 
we cannot expect strong convergences (as in part $(iii)$ of the theorem) unless~$\operatorname{supp}(\rm Q)$ is strongly compact.\end{itemize}

Theorem~\ref{Theorem stability} also has important consequences and  potential applications.
\begin{itemize}
    \item[(d)] [Discrete ${\rm P}_n$] When ${\rm P}_n$ is discrete,   a transport map inducing the coupling $\gamma_n\in \Pi({\rm P}_n,{\rm Q}_n)$  may not exist, but  (under second order moment assumptions) the solution $\gamma_n$  of the optimal transport problem $$ 
\inf_{\pi\in \Pi( {\rm P}_n,{\rm Q}_n)} \int \|x-y \|^2d\pi(x,y)
 $$  always exists. In the notation of Theorem~\ref{Theorem stability},  the sequence $\gamma_n$ of couplings  thus still provides a consistent estimator of $\nabla\psi$ since, for any $h \in \H$,
$$
 \sup_{{ (x,y)\in \,{\rm supp}(\gamma_n), x\in K} } \langle y-\nabla\psi(x), h \rangle  \longrightarrow 0 \quad\text{ as $n\to\infty$}.
$$
  \item[(e)]  [Glivenko-Cantelli] When ${\rm P}_n$ and ${\rm Q}_n$ are supported on two disjoint sets   of the same cardinality and give the same probability mass to each point (i.e., ${\rm P}_n$ and ${\rm Q}_n$ are the empirical measures over these two sets), then, for each $n\in \N$ the transport problem~\eqref{OTIntro} from ${\rm P}_n$ to ${\rm Q}_n$  admits a solution  $T_n$ defined uniquely on the support points of ${\rm P}_n$ only.  In this case, under the  assumptions of Theorem~\ref{Theorem stability}, for any $h \in \H$,
\begin{equation}\label{eq:GC}
\sup_{{ x\in \,{\rm supp}({\rm P}_n) \cap K} } \langle T_n(x)-\nabla\psi(x), h \rangle  \longrightarrow 0\quad\text{ as $n\to\infty$}.
\end{equation}
This scenario frequently arises  in statistical applications, where ${\rm P}_n$ is the empirical measure associated with an i.i.d.~sample $X_1, \dots, X_n\sim {\rm P}$ and ${\rm Q}_n$ is a certain discretization of a given reference measure ${\rm Q}$. 
For $\H = \R$ and ${\rm Q} = $ Uniform[0,1], the transport map $T_n$ then reduces to the usual cumulative distribution function:   \eqref{eq:GC}, therefore, has the form of  an extended local Glivenko-Cantelli  theorem (where the convergence, moreover, holds $\rm P$-a.s.). To obtain a ``full''  (non-local, i.e., over all of $\H$) Glivenko-Cantelli result, we may investigate the regularity properties of $\nabla\psi$ under some smoothness assumptions on ${\rm P}$. In the Euclidean case, such results have been established in~\cite{Hallin2020DistributionAQ},~\cite{DELBARRIO2020104671},~\cite{GhosalSenAOS} and \cite{FIGALLI2018413}. Their proof uses  the well-known Caffarelli theory (see \cite{FigalliBook} and references therein) which, however, has not been fully developed in the infinite-dimensional case.
\end{itemize} 

\subsection{A central limit result for Wasserstein distances}\label{CLTsec}
The stability \eqref{(ii)} of  potentials implies (under appropriate moment assumptions), via the Efron-Stein-inequality-based argument of \cite{delBarrioLoubes19} and \cite{delbarrio2021central}, a central limit behavior of the fluctuations of the squared  Wasserstein distance  
$$ 
 \mathcal{T}_2({\rm P}_n,{\rm Q}) \coloneqq  \inf_{\pi\in \Pi( {\rm P}_n,{\rm Q})} \int \|x-y \|^2d\pi(x,y),
 $$ 
between the empirical distribution of   an i.i.d.\ sample $X_1, \dots, X_n$ from $\rm P$  and $\rm Q$.  In the following theorem, we show that, under suitable conditions, the fluctuations of $\mathcal{T}_2({\rm P}_n,{\rm Q})$ around its expectation are asymptotically Gaussian.

\begin{Theorem}\label{TCL}
Let $\rm P, \rm Q\in\P(\H)$ be such that ${\rm P}\in \Pl(\H)$ 
and has connected support. Assume furthermore that the boundary $\partial\, {\operatorname{supp}({\rm P})}$ of ${\operatorname{supp}({\rm P})}$  has  $\rm P$-probability zero, 
that~$\int \|x\|^4 d {\rm P}(x) <~\!\infty$, and that~$\operatorname{supp}({\rm Q})$ is bounded. Then, there exists a unique (up to additive constants) proper l.s.c.\ convex function $\psi$  such that   $(\nabla \psi )\# \rm P={\rm Q}$, and
\begin{align}\label{conclusion1}
\sqrt{n}\Big(\mathcal{T}_2({\rm P}_n, {\rm Q})-{\E} [\mathcal{T}_2({\rm P}_n, {\rm Q})]\Big)\stackrel{w}{\longrightarrow} N\Big(0, \sigma^2_2(\rm P,\rm Q)\Big), 
\end{align}
where $N(0, \sigma^2_2(\rm P, \rm Q))$ is a univariate Gaussian distribution with mean zero and variance
\begin{align}\label{sigma_def}
\sigma^2_2({\rm P},{\rm Q})\coloneqq\int{(\|x\|^2-2\psi({x}))^2}{\rm dP}({x})-\left( \int{\left(\|x\|^2-2\psi({x})\right)}{\rm dP}({x})\right)^2.
\end{align}
\end{Theorem}

Note that a central limit theorem for $ \mathcal{T}_2({\rm P}_n,\rm Q)$ centered at its population counterpart $ \mathcal{T}_2(\rm P,\rm Q)$ is, in general, impossible in infinite-dimensional Hilbert spaces due to the well-known curse of  dimensionality: see, e.g.,~\cite{weedBach}, who show that the  convergence to zero of the bias is much  slower than the decrease of the variance in dimension higher than $5$.\footnote{There are a few  exceptions, though. In  some particular cases, when either $\rm P$ or $\rm Q$ or both are supported on simpler spaces (such as a finite number of points or a lower-dimensional  manifold), replacing the expectation of $ \mathcal{T}_2({\rm P}_n,{\rm Q})$ with the population quantity $ \mathcal{T}_2({\rm P},{\rm Q})$ in~\eqref{conclusion1} is possible; see the results of~\cite{delbarrio2021centraldisc} and \cite{Hundrieser2022EmpiricalOT}.}

The proof of Theorem~\ref{TCL} mainly consists of  showing the  $\rm P$-a.s.~uniqueness of $\psi$;  the rest of the proof follows, almost \emph{verbatim}, along the arguments developed in  \cite{delbarrio2021central}. Details, thus, are  skipped. To prove the $\rm P$-a.s.\ uniqueness of $\psi$, we first show that, if the subdifferentials of two convex functions coincide on a dense subset of {an open convex set},   they coincide  everywhere on that  {open set}. The  following lemma makes this precise and is inspired by the proof of  Case (2-c) in \cite{CorderoErausquin2019RegularityOM}. 
\begin{figure}[h]$\,$\vspace{4mm}
    \centering
    \includegraphics[width=0.95\textwidth]{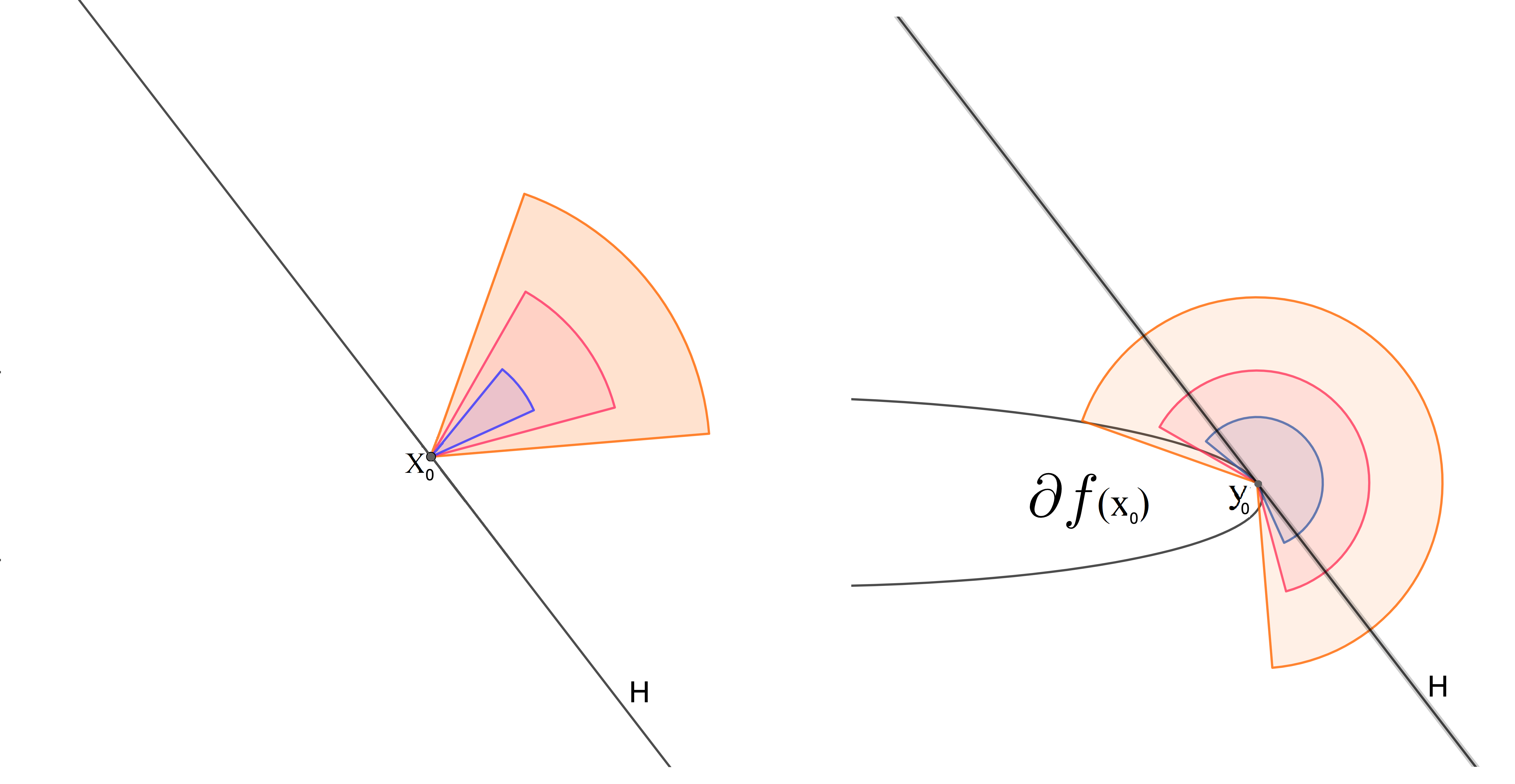}
    \caption{A visual illustration of the proof of Lemma \ref{Lemma:uniqueness}. The subdifferential $\partial f$ maps each region of the figure on the left to the   region of the same color on the right. In each region we can pick  a point in~$ \mathcal{D}_{f\cap g}\coloneqq  \{ x \in \H:\  \partial f(x)\cap  \partial g(x)\neq \emptyset \} $,  hence create a sequence $x_n\to x_0$ with~$\{x_n\}_{n\in \N}\subset\mathcal{D}_{f\cap g}$. Then  the only possible limit for $\{y_n\}_{n\in \N}$ where $y_n\in  \partial f(x_n)\cap  \partial g(x_n)$ is~$y_0$. }
    \label{fig}
\end{figure}
\begin{Lemma}\label{Lemma:uniqueness}
Let $f, g:\H\to (-\infty,+\infty]$ be proper l.s.c. convex functions such that the set
$$\mathcal{D}_{f\cap g}\coloneqq \{ x \in \H:\  \partial f(x)\cap  \partial g(x)\neq \emptyset \} $$
is dense in an open convex set $\mathcal{U}\subseteq\operatorname{int}(\operatorname{dom}(f))\cap \operatorname{int}(\operatorname{dom}(g))$.  Then $ \partial f=\partial g $ in $\mathcal{U}$. 
\end{Lemma}
\begin{proof}
Suppose that there exists $x_0\in \mathcal{U}$ such that $ \partial f(x_0)\neq \partial g(x_0) $. We will show that this assumption leads at a contradiction. The set $\partial f(x_0)$ is convex and weakly compact~\citep[Proposition 16.14]{bauschkeMonotoneHilbert} so that, via the Lindenstrauss theorem (see~\citet[Theorem 4]{Lindenstrauss1963OnOW}),  $\partial f(x_0)$ is the closure of the convex hull of the set~$S(x_0)$ of its exposed points. That is,  for any $y_0\in S(x_0)$, there exists a supporting hyper
plane~$H\coloneqq\{z \in \H:\langle z,b \rangle+a=0\}$ such that $$\partial f(x_0)\cap H = \{y_0\} \quad \mbox{while} \quad \partial f(x_0)\subseteq H^-,$$ where~$H^-\coloneqq \{z \in \H:\langle z,b \rangle+a\leq 0\}$; see the visual illustration  provided in Figure~\ref{fig}.  This implies that 
\begin{equation}\label{H+}
    \partial f(x_0)\cap H^+= \partial f(x_0)\cap \left((\H\setminus H^-)\cup H \right) =  \{y_0\}, 
\end{equation}
where $H^+\coloneqq \{z \in \H:\langle z,b \rangle+a\geq  0\}$.

Without loss of generality, we can assume~$x_0=y_0=0$, $\|b\|=1$, and $a=0$ (this can be achieved by redefining $f$ as in~\eqref{eq:g*}).
Denote by $\pi$ be the orthogonal projection to the space~$\{ \lambda \, b: \ \lambda\in \R \} ^{\bot}$. Thus, for any $z \in \H$, we can write $z = \langle z,b\rangle b + \pi(z)$ where $\langle b, \pi(z) \rangle = 0$. Consider the open convex truncated cone (see Figure~\ref{fig} for an illustration)
$$C_{n}\coloneqq \left\{z \in \H: {n}^{-1}  \langle z,b \rangle> \|\pi(z) \|\right\}\cap  \mathcal{B}(0,\ {n}^{-1})\cap \mathcal{U}.$$  
Let  $\{x_n \}_{n \ge 1}\subseteq \operatorname{dom}(f)\cap \operatorname{dom}(g) $ be a sequence such that $x_n\in  C_{n} $ for all $n\in \N$.\linebreak Clearly, $x_n\to 0 \equiv x_0$. However, by strong-to-weak  continuity of the subdifferential (see part $(i)$ of  Lemma~\ref{Lemma:strongtoweak}), any sequence $y_n\in \partial f(x_n)\cap \partial g(x_n)$ converges  weakly (possibly along  subsequences) to\linebreak some~$u\in \partial f(x_0) \cap \partial g(x_0)$. Let us now show that $u=y_0=0$. 
To do so,   observe that, by cyclical monotonicity, 
$$0\le\langle x_n - x_0, y_n -y_0 \rangle = \langle x_n, y_n \rangle \langle \pi(x_n), \pi(y_n)\rangle+\langle x_n, b\rangle\langle b, y_n\rangle \leq   \| \pi(x_n)\|\| \pi(y_n)\|+\langle x_n, b\rangle\langle b, y_n\rangle .$$ 
Since $x_n\in  C_{n} $,  $$ \langle x_n, b\rangle ({n}^{-1} \| \pi(y_n)\|+\langle b, y_n\rangle)> \| \pi(y_n)\|\|\pi(z) \|+\langle x_n, b\rangle \langle b, y_n\rangle \geq 0$$
with $ \langle x_n, b\rangle> 0$, so that 
$ y_n\in\{y \in \H: \ -{n}^{-1} \|\pi(y)\|\leq \langle y,b\rangle \} .$ 

Let us show that the weak limit~$u$ of~$\{y_n\}_{n\in \N}$ is such that  $\langle u,b \rangle \ge 0$. 
Suppose therefore  that~$\langle u,b \rangle =-\delta<0$. Then, by the definition of the weak limit of $\{y_n\}$, we obtain~$\langle y_n,b \rangle \to -\delta$; also note that $\lim\sup\|\pi(y_n)\|\leq \lim\sup \|y_n\|<\infty$ (see e.g.,~\citet[Proposition 3.13 (iii)]{Brezis}).  We also know that $-{n}^{-1}\|\pi(y_n)\|\leq \langle y_n,b\rangle$, where   taking limits yields the contradiction
$ 0\leq -\delta$. Hence, $\langle u,b \rangle \ge 0$.

As a consequence, $u\in \partial f(x_0)\cap  \partial g(x_0)\cap H^+$ where $H^+ \coloneqq  \{z:\langle z,b \rangle \geq  0\}$. However, in view of \eqref{H+},   $y_0$ is the only point in~$ \partial f(x_0)\cap H^+ $, so that $y_0=u$. This means that any $y_0$ belonging to $S(x_0)$ (the set extreme points  of $\partial f(x_0)$)  also belongs to $\partial g(x_0)$. Hence, $S(x_0)\subseteq \partial g(x_0)$. Since~$\partial g(x_0)$ is convex and weakly compact, $\partial f(x_0)\subseteq \partial g(x_0)$. The converse follows by noting that $f$ and $g$ are playing fully symmetric roles. Thus we obtain that $\partial f(x_0)= \partial g(x_0)$, which yields the contradiction. The claim follows.
\end{proof}

\begin{proof}[Proof of Theorem~\ref{TCL}] 
The proof of Theorem~\ref{TCL} is  similar to that of~\citet[Theorem~4.5]{delbarrio2021central} (see also \cite{gonzalezdelgado2021twosample} and \cite{delBarrioLoubes19})  details therefore are skipped and only a brief outline of the proof is given. 

The proof proceeds via an {application of the Efron-Stein inequality to the random variable $$R_n \coloneqq \mathcal{T}_2({\rm P}_n, {\rm Q}) - \int \psi(x) d{\rm P}_n(x).$$}  
Defining  independent copies $X'_1, \ldots, X'_n$ of the random variables $X_1, \ldots, X_n$,  denote by ${\rm P}_n^{(i)}$ the empirical measure on $(X_1, \ldots, X_{i-1}, X'_i, X_{i+1}, \ldots, X_n)$ and by~$R_n^{(i)}$   the version of $R_n$ computed from ${\rm P}_n^{(i)}$ instead of ${\rm P}_n$.  The proof consists in showing that the sequence $n(R_n - R_n^{(i)})$  converges almost surely  to $0$ while  $n^2 \E[(R_n - R_n^{(i)})^2]$ is bounded by some constant $M$. This latter claim is a consequence of  the finite fourth-order moment assumptions on~${\rm P}$ and~${\rm Q}$ while the first one (see Eq.~(6.2) in \cite{delbarrio2021central}) results from the stability  of the potentials in Theorem~\ref{Theorem stability}~$(ii)$---provided, however, that it holds true in this context. Assuming that it does,  it follows from the Banach-Alaoglu theorem (see e.g.,~\citet[Theorem 3.16]{Brezis}) and the Banach-Saks property (see e.g.,~\citet[Exercise 5.24]{Brezis})  that there exists a subsequence of~$\{n|R_n - R_n^{(i)}|\}_{n \in \mathbb{N}}$ the Ces\`{a}ro mean\footnote{T
Recall that the Ces\`{a}ro mean of a sequence $\{ x_n\}_{n\in \N}$ is the sequence $\{ y_n\coloneqq \frac{1}{n}\sum_{i=1}^n x_i\}_{n\in \N}$.} of which converges to $0$ in $L^2(\mathbb{P})$. The same property holds for the Ces\`{a}ro means of  subsequences of $\{\sqrt{n}(R_n - \E R_n)\}_{n \in \mathbb{N}}$.  This implies (see p.~21 in \cite{delbarrio2021central}) that~$\sqrt{n}\E|R_n - \E R_n|$ converges to zero and the central limit theorem holds.

 The only ingredient missing in this sequence of arguments, thus, is the stability result of the potentials (Theorem~\ref{Theorem stability}-$(ii)$). The   requirement here is  the uniqueness (up to additive constants) of the population potential $\psi$, which we now establish: given two proper l.s.c.\  convex functions $\psi$ and~$f$ such that $ \nabla \psi=\nabla f $ on a dense subset of $\operatorname{int}(\operatorname{supp}(\rm P)) \subset \H$, let us  show that $\psi$ and $f$ are equal (up to an additive constant) on $\operatorname{int}(\operatorname{supp}(\rm P))$.

Let $p\in \operatorname{int}(\operatorname{supp}(\rm P))$. Then there exists an open convex neighborhood $\mathcal{B} \subset \operatorname{int}(\operatorname{supp}(\rm P))$ of~$p$.  We have that $ \nabla \psi=\nabla f $ on a dense subset~$\mathcal{D}$ of $\mathcal{B}$. Lemma~\ref{Lemma:uniqueness} thus  implies that  $\partial {\psi}=\partial {f} $ on~$\mathcal{B}$. Let $\mathcal{B}'$ be an open convex neighborhood of $p$ such that its topological closure $\overline{\mathcal{B}'}$ is contained in~$\mathcal{B}$. Define $s_{\overline{\mathcal{B}'}}$ as the support function of $\overline{\mathcal{B}'}$, i.e., $s_{\overline{\mathcal{B}'}}(x)=0$ if $x\in \overline{\mathcal{B}'}$ and $+\infty$ otherwise.  Since~$ \partial \tilde{\psi}(h)=\partial \tilde{f}(h) = \emptyset$ for $h\in \H\setminus \overline{\mathcal{B}'}$  and $ \partial \tilde{\psi}(h)=\partial \tilde{f}(h)$ for $h\in \overline{\mathcal{B}'}$,  the functions~$\tilde{\psi}\coloneqq \psi +s_{\overline{\mathcal{B}'}}$ and $\tilde{f}\coloneqq f +s_{\overline{\mathcal{B}'}}$ are proper l.s.c.\ convex functions with $\partial \tilde{\psi}=\partial \tilde{f} $ on $\H$ and   $ \nabla \tilde{\psi}=\nabla \tilde{ f} $  in the dense subset $\tilde{\mathcal{D}} \coloneqq \mathcal{D} \cup (\H \setminus \overline{\mathcal{B}})$ of $\H$. 
\citet[Theorem B]{RockafellarMaximalMonot} then yields the existence of some~$a=a_{p}\in 
\R$ such that  $\tilde{\psi} =\tilde{ f}+a$. Thus, $\psi=f+a$ on $\mathcal{B}'$. 

Up to this point, we have proven that, for all $p\in \H$, there exists a neighborhood $\mathcal{B}'$ of $p$ and a constant $a=a_p\in \R$ such that $\psi=f+a$ on $\mathcal{B}'$.  Using a connectedness\footnote{Recall that a set is connected if it cannot be written as the union of two non-empty open and disjoint sets.} argument, let us show that  this constant  actually does not depend on $p$.  The set~$\Theta$ of all $x$ in $\operatorname{int}(\operatorname{supp}(\rm P))$ such that~$\psi(x)=f(x)+a$ {is open (i.e.,  each $q\in \Theta$ has a neighborhood where  $\psi=f+a$) and non-empty (since $p\in \Theta$). Its complement is open, too (for each $q\not\in \Theta$ there exists a neighborhood of $q$ where~$\psi=f+b$, with $b\neq a$, and  this neighborhood obviously does not contain any element of~$\Theta$). Since $\operatorname{int}(\operatorname{supp}(\rm P))$ is connected, $\Theta=\operatorname{int}(\operatorname{supp}(\rm P))$. This completes the proof of Theorem~\ref{TCL}.} 
\end{proof}

\subsection{Proof of Theorem \ref{Theorem stability} }\label{Theorem stabilitySec}
The proof of Theorem~\ref{Theorem stability} relies on a series of lemmas. {Some of these lemmas are self-contained ``general'' results; some others address the specific setting of Theorem~\ref{Theorem stability}.}

\begin{Lemma}\label{lem:Coupling}
    Under the assumptions of Theorem \ref{Theorem stability}, $\gamma_n\xrightarrow{w}\gamma$, where $\gamma\in \Pi({\rm P},{\rm Q})$  is the unique probability measure  such that 
 $\operatorname{supp}(\gamma)\subseteq \partial \psi  $  for some l.s.c.\ convex function $\psi$. 
\end{Lemma}
\begin{proof}Lemma~\ref{Lemma:MCann} implies that for any subsequence $\{\gamma_{n_k} \}_{k\in \N}$ there exists a further sub\-se-\linebreak quence~$\{\gamma_{n_{k_i}} \}_{i\in \N}$ converging weakly  to a probability measure $\gamma\in \Pi({\rm P},{\rm Q})$  such that 
 $\operatorname{supp}(\gamma)\subseteq \partial \psi  $  for some l.s.c.\ convex function $\psi$. Since  $\partial \psi$ is ${\rm P}$-a.s.~a singleton, it follows from Theorem~\ref{Theorem:Maccan} that~$\gamma = ({\rm Identity}, \nabla \psi)\# {\rm P}$ is unique. Hence $\gamma_n\xrightarrow{w}\gamma$.
\end{proof}

In order to move from weak convergence of   couplings to convergence of mappings, we use the set-topology of  subdifferentials of the form $\partial\psi_{n}$ where $\psi_{n}$ denotes a sequence of convex functions defined over $\cal H$. Denoting by $\{A_n\}_{n\in \N}$  a sequence of subsets of a  {\it second countable topological space}\footnote{{A topological space is said to be {\it second countable} if its topology admits   a countable basis.}} $\mathcal{Y}$,  define the {\it inner} and {\it outer limits} of $\{A_n\}_{n\in \N}$ as
\begin{align*}
& \operatorname{Liminn}_{n}^{\mathcal{Y}}A_{n}\coloneqq \{ x\in\mathcal{Y}: \ \ \text{exists  $\{x_n\}_{n\in\N}$ with $x_n\in A_{n}$ \text{such that} $x_n\xrightarrow{\mathcal{Y}}x$}\text{ as $n\to\infty$}\}\vspace{-4mm}
\end{align*}
and\vspace{-3mm}
\begin{align*}
&\hspace{7.5mm}\operatorname{Limout}_{n}^{\mathcal{Y}}A_{n}\coloneqq \{ x\in\mathcal{Y}: \ \ \text{exists  $\{x_{n_k}\}_{k\in\N}$ with $x_{n_k}\in A_{n_k}$ \text{such that} $x_{n_k}\xrightarrow{\mathcal{Y}} x$}\text{ as $k\to\infty$}\},
\end{align*}
respectively. Here the convergence $x_n\xrightarrow{\mathcal{Y}}x$ is to be understood in the sense of the topology of $\mathcal{Y}$, i.e., for any neighborhood $\mathcal{U}_x$ of $x$, there exists $N_{\mathcal{U}_x}\in \N$ such that $x_n\in \mathcal{U}_x$, for all $n\geq N_{\mathcal{U}_x}$.   The corresponding  {\it limit} exists if and only if the inner  and  outer limits coincide, in which case we say   that~$\{A_n\}_{n\in \N}$ {\it converges in the  Painlev\'{e}-Kuratowski sense} with respect to the topology of $\mathcal{Y}$.  When  the space  $\mathcal{Y}$ is not second countable, the inner and outer limits are usually expressed in terms of nets instead of sequences. For ease of reference, we recall here a fundamental result on this type of set convergence, known in the literature as Mr\'{o}wka's theorem (see, e.g.,~\citet[{Theorem 5.2.12}]{Beer1993TopologiesOC}).  
\begin{Lemma}[Mr\'{o}wka]\label{Mrowka}
    Let $\mathcal{Y}$ be a  second countable topological space. Then, for any sequence  $\{A_n\}_{n\in \N}
    $ of subsets of $\mathcal{Y}$,{ there exists a Painlev\'{e}-Kuratowski
convergent subsequence $\{A_{n_k}\}_{k\in \N}$}. 
\end{Lemma}

In particular, for $\mathcal{Y}=\H\times \H$,  denote by
\begin{align*}
&\operatorname{Limout}_{n}^{s-w} A_{n}\coloneqq \{ (x,y)\in \H\times \H: \,\text{exists  $(x_{n_k}, y_{n_k})\in A_{n_k}$ \text{such that} $x_{n_k}\to x$ and $y_{n_k}\rightharpoonup y$}\}\vspace{-4mm}
\end{align*}
and\vspace{-2mm}
\begin{align*}
\hspace{-2.5mm}&
    \operatorname{Liminn}_{n}^{s-w} A_{n}\coloneqq \{ (x,y)\in \H\times \H: \ \ \text{exists  $(x_n, y_n)\in A_{n}$ \text{such that} $x_n\to x$ and $y_n\rightharpoonup y$}\},
\end{align*}
 respectively, the {\it strong-to-weak outer limits} and {\it strong-to-weak inner limits} of $\{A_n\}_{n\in \N}$. For sequences of subdifferentials of the form~$\partial\psi_n$, when $\operatorname{Liminn}_{n}^{s-w}\partial\psi_n=\operatorname{Limout}_{n}^{s-w}\partial\psi_n$,  we denote by $\operatorname{Lim}_{n}^{s-w}\partial\psi_n$  the strong-to-weak  {\it Painlev\'{e}-Kuratowski limit}.
Painlev\'{e}-Kuratowski limits are appropriate for sequences of cyclically monotone sets since, as the following result shows, they preserve that property. 
\begin{Lemma}\label{Lemma:LimitCyclically}
Let $\{\Gamma_n\}_{n\in \N}$  
 be a sequence of cyclically monotone subsets of the product space~$\H\times \H$. If $\Gamma^{s-w}\coloneqq \operatorname{Liminn}_{n}^{s-w}\Gamma_n\neq \emptyset$, then $\Gamma^{s-w}$ is also cyclically monotone. 
\end{Lemma}
\begin{proof}
By definition, for each finite $N$-tuple $\{(x_k,y_k)\}_{k=1}^N\subseteq \Gamma^{s-w}$, there exists a sequence\linebreak of~$N$-tuples  $\{\{(x_k^n,y_k^n)\}_{k=1}^N\}_{n\in \N}$ with  $\{(x_k^n,y_k^n)\}_{k=1}^N\subseteq \Gamma_{n}$ such that 
$x_k^n\rightarrow x_k$ and $y_k^n\rightharpoonup y_k$ as~$n \to \infty$, for all $k=1, \dots, N$. For each $n\in \N$, cyclical monotonicity of $\Gamma_n$ implies 
\begin{equation}
    \label{takingLim}
    	\sum_{k=1}^N \langle x_{k}^n,y_{k+1}^n-y_k^n\rangle\leq 0.
\end{equation}
Now, { $\{\|y_{k+1}^n-y_k^n\|\}_{n\in{\mathbb N}}$ is bounded (e.g.,~\citet[Proposition 3.13 (iii)]{Brezis})}, so that
\begin{multline*}
    |\langle x_{k}^n,y_{k+1}^n-y_k^n\rangle-\langle x_{k},y_{k+1}-y_k\rangle|= |\langle x_{k}^n-x_{k},y_{k+1}^n-y_k^n\rangle +\langle x_{k},y_{k+1}^n-y_k^n-(y_{k+1}-y_k)\rangle|\\
    \leq \| x_{k}^n-x_{k}\|\|y_{k+1}^n-y_k^n\| +|\langle x_{k},y_{k+1}^n-y_k^n-(y_{k+1}-y_k)\rangle|\to 0.
\end{multline*}
Taking limits in \eqref{takingLim} yields 
$
	\sum_{k=1}^N \langle x_{k},y_{k+1}-y_k\rangle\leq 0
$. 
The claim follows.
\end{proof}

Let us now get back to the setting and the notation of Theorem~\ref{Theorem stability}. Observe that $\cal H \times \cal H$ with strong topologies on both sets is a separable metric space and hence is second countable (see~\citet[Exercise 2.23]{rudinMathematicalAnal}). The space $\H$ endowed with the weak topology, however, is not second countable  (see \citet[Corollary 1]{Helmberg}), so that Mr\'{o}wka's theorem does not directly apply to $ \operatorname{Lim}_{n}^{s-w}\partial \psi_n $ in the product space~$ \H\times \H$   when considering the strong topology in the first factor  and the weak topology in the second one---see \citet[Proposition 5.2.13]{Beer1993TopologiesOC} for a counter-example assuming the continuum hypothesis. However, consider the open ball $\mathcal{B}(0,{M})$ with radius $M$ centered at $0$, endowed with  the metric structure~$(\mathcal{B}(0,{M}), \| \cdot\|_w)$, 
where
\begin{equation}
    \label{metric1}
    \|x \|_{w}\coloneqq\sum_{k\in \N} \frac{1}{2^k}\frac{| \langle x, e_k \rangle |}{1+| \langle x, e_k \rangle |}
\end{equation}
with $\{e_k \}_{k\in \N}$  an orthonormal basis of $\cal H$. The norm~\eqref{metric1} metrizes weak convergence in~$\mathcal{B}(0,M) $, so that if $\bigcup_{n\in \N, x\in K}\partial \psi_n(x)$ is a subset of $\mathcal{B}(0,M)$ --- which holds by assumption --- it is relatively compact in $(\mathcal{B}(0,M), \| \cdot\|_{w})$. Moreover, $(\mathcal{B}(0,M), \| \cdot\|_{w})$ is separable and, being metrizable, it is second countable, see \citet[Exercise 2.23]{rudinMathematicalAnal}. Mr\'{o}wka's theorem, thus, now  applies. 
   
Mr\'{o}wka's theorem (Lemma~\ref{Mrowka})  implies that, given the sequence~$\{\partial\psi_{n}\}_{n\in \N}$, there exists a   subsequence $\{\partial\psi_{n_k}\}_{k\in \N}$ which converges, in the  Painlev\'{e}-Kuratowski sense with respect to the strong-to-weak topology, to some $ \Gamma^\prime $. But this limit may be the empty set, rendering the application of Mr\'{o}wka's theorem meaningless. The following Lemma shows that this is not the case for the  sequence of subdifferentials. The strong-to-weak  Painlev\'{e}-Kuratowski inner limit  $\Gamma$ of~$\{\partial\psi_{n}\}_{n\in \N}$ --- hence, also the set~$\Gamma^\prime\supseteq\Gamma$ --- is non-empty and contains all the points of $\operatorname{supp}(\rm P)$ at which $\partial \psi$ is a singleton.  
\begin{Lemma}\label{Lemma:SetLim}
Under the assumptions of Theorem \ref{Theorem stability},   
the set $\operatorname{Liminn}_{s-w}\partial \psi_n$ in nonempty and, moreover, contains~$\{  (h, \nabla \psi(h)):\ h\in \operatorname{supp}(P)\cap \operatorname{dom}(\nabla \psi)\}$. 

\end{Lemma}
\begin{proof}
Let $h\in \operatorname{supp}(P)\cap \operatorname{dom}(\nabla \psi) $ and let $y \coloneqq \nabla \psi(h)$. Denote by $\{\mathcal{U}_{y}^m\}_{m\in \N}$  a decreasing (i.e.,~$\mathcal{U}_{y}^m\subseteq \mathcal{U}_{y}^{m-1}$) countable sequence of   neighborhoods (with respect to the strong topology) of~$y$ such that $\bigcap_{m\in \N} \mathcal{U}_{y}^m=\{ y\}$. Then,~\citet[Theorem 21.22]{bauschkeMonotoneHilbert} entails the existence of a {\it selection}\footnote{A {\it selection} of $\partial\psi$ is a map $\mathcal{Q}:  \operatorname{dom}(\partial \psi) \to \H$ such that $(x,\mathcal{Q}(x))\in \partial \psi$ for all $x\in \operatorname{dom}(\partial \psi)$ (see \citet[p.~2]{bauschkeMonotoneHilbert}).}~$\mathcal{Q}$ of $\partial\psi$ and a decreasing sequence of  neighborhoods (with respect to the weak topology) $\{\mathcal{V}_{h}^m\}_{m\in \N}$ of $h$ such that  $\bigcap_{m\in \N} \mathcal{V}_{h}^m=\{ h\}$ and~$\mathcal{Q}(\mathcal{V}_{h}^m)\subseteq \mathcal{U}_{y}^m $  for all $m\in \N$. Since~$\mathcal{V}_{h}^m$ has non-empty interior and $h \in \operatorname{supp}(P)$, and since~$\gamma=({\rm Identity} \times \mathcal{Q} ){\#}{\rm P}$, we have 
$$  \gamma(\mathcal{V}_{h}^m\times   \mathcal{U}_{y}^m)\geq \gamma(\mathcal{V}_{h}^m\times  \mathcal{Q}(\mathcal{V}_{h}^m))=\mu(\mathcal{V}_{h}^m)=:\delta_m>0.$$ 
Moreover, $\delta_m\to 0$ monotonically  because $\bigcap_{m\in \N} \mathcal{V}_{h}^m=\{ h\}$ and $\mathcal{V}_{h}^m\subseteq \mathcal{V}_{h}^{m-1}$. Fix  $m\in \N$:  since~$\gamma_n$ converges weakly to  $\gamma=({\rm Identity} \times \mathcal{Q} ){\#}{\rm P}$, the Portmanteau theorem yields 
$$ \lim\inf_n\gamma_n(\mathcal{V}_{h}^m\times \mathcal{U}_{y}^m)\geq \delta_m, $$
so that there exists $ N_m\in \N$ and a sequence $(h_{n},y_{n})$ such that $$(h_{n},y_{n})\in (\mathcal{V}_{h}^m\times  \mathcal{U}_{y}^m)\cap \operatorname{supp}(\gamma_n)$$  for all $n\geq N_m$. By definition of $\gamma_n$, there exists  $(h_n,y_{n})\in \partial \psi_n$ such that $h_n\in \mathcal{V}_{h}^m$ and~$y_n\in~\!\mathcal{U}_{y}^m$. Letting~$m\to~\!\infty$ yields a sequence 
   $\{(h_{n}, y_n)\}_{n\in \N}$ with $(h_{n}, y_n)\in  \partial \psi_n $, $h_{n}\longrightarrow h$, and~$ y_n\rightharpoonup~\!\nabla\psi(h)$, which completes the proof.\end{proof}

Lemma~\ref{Lemma:LimitCyclically} implies that, in the context of Theorem~\ref{Theorem stability},  for any subsequence $ \{ n_k\}_{k\in \N}$, the limiting sets  {$ \Gamma\coloneqq \operatorname{Liminn}_{n}^{s-w}\partial\psi_{n}$ and  $\Gamma^\prime\coloneqq~\!\operatorname{Lim}_{k}^{s-w}\partial\psi_{n_k}$} both  are cyclically monotone. As a consequence, there exists a convex function $\rho$ such that $\Gamma\subseteq \Gamma^\prime\subseteq \partial \rho$ (see \citet[Theorem~B]{RockafellarMaximalMonot}). In view of Lemma~\ref{Lemma:SetLim}, we have, still in the setting of Theorem~\ref{Theorem stability},~$(h, \nabla\psi(h))\in~\!\partial \rho$ for any $h\in \operatorname{supp}(P)\cap \operatorname{dom}(\nabla \psi)$. 
 Since this set is dense in $\operatorname{supp}({\rm P})$ (see \citet[Theorem 21.22]{bauschkeMonotoneHilbert}), Lemma~\ref{Lemma:uniqueness} entails $\partial \rho=\partial \psi$ in $\operatorname{int}\operatorname{supp}(\rm P)$, so   that~$\Gamma\subseteq \Gamma^\prime\subseteq \partial\psi$. This constitutes a fundamental difference with   the finite-dimensional case. In Euclidean spaces, indeed, the limits of maximal monotone operators (subdifferentials of proper l.s.c.\ convex functions) is automatically maximal monotone (see e.g.\  \cite{adlyRockafellar}) and, instead of~$\Gamma\subseteq \Gamma^\prime\subseteq \partial\psi$, it holds that  $\Gamma= \Gamma^\prime= \partial\psi$.  In the infinite-dimensional case, with the notation of Theorem~\ref{Theorem stability}, we thus have the following property.
 
\begin{Lemma}\label{Lemma:Subdiff}
    Under the assumptions of Theorem \ref{Theorem stability},  
    \begin{enumerate}[(i)]
    \item for any  $x\in \operatorname{dom}(\nabla \psi) \cap\, \operatorname{int}(\operatorname{supp}(\rm P))$, there exists a sequence $\{(x_n,y_n)\}_{n \ge 1}$ such that~$x_n\to x$ and $y_n\rightharpoonup \nabla \psi(x)$ with $y_n\in \partial \psi(x_n)$ for $n$ large  enough;
     \item for any sequence $\{(x_n,y_n)\}$ (or subsequence thereof) such that $x_n\to x\in \operatorname{int}(\operatorname{supp}(\rm P))$ and~$y_n\rightharpoonup y$, with~$y_n\in \partial \psi(x_n)$ for $n$ large  enough,   $ y\in \partial\psi(x) $.  
\end{enumerate}
\end{Lemma}

We are now ready for the proof of Theorem~\ref{Theorem stability}. 

\begin{proof}[{\bf Proof of  Theorem~\ref{Theorem stability}}]
 First consider parts $(i)$ and $(iii)$. Set $h\in \H$ and suppose that 
$$ \liminf_{n\to \infty}\sup_{(x,y)\in \partial \psi_n, \, x \in K} |\langle y-\nabla \psi(x), h \rangle | > \epsilon$$  
for some $K\subseteq\operatorname{dom}(\nabla \psi)\cap \operatorname{int}(\operatorname{supp}(P))$: that implies the existence of sequences $\{x_{n_k} \}_{k\in \N}\subseteq K$ and $\{y_{n_k} \}_{k\in \N}\subseteq \H$ with $y_{n_k}\in \partial \psi_{n_k}(x_{n_k})$ such that, for all $k \ge 1$,
\begin{equation}
    \label{ToContradice}
    |\langle y_{n_k}-\nabla \psi(x_{n_k}), h \rangle | >\epsilon.
\end{equation}
Since (by assumption) $\bigcup_{n\geq N_0, x\in K}\partial \psi_n(x)$  is contained in the ball $\mathcal{B}(0,M)$, there exist (by the Banach-Alaoglu theorem; see e.g.,~\citet[Theorem 3.16]{Brezis}) a weak limit~$ y$ of the sub\-sequence~$\{y_{n_{k_i}} \}_{i\in \N}$  and (by the strong compactness of $K$) a strong limit $x$ of the subse-\linebreak quence~$\{x_{n_{k_i}} \}_{k\in \N}$. The limit $(x,y)$ belongs to the set $\operatorname{Liminn}_{i}^{s-w} \partial \psi_{n_{k_i}}$, hence, via Lemma~\ref{Lemma:Subdiff},   the only possible value of $(x,y)$ is $(x,\nabla\psi(x))$. This yields a contradiction to~\eqref{ToContradice} and proves part~$(i)$ of the theorem, of which part~$(iii)$ is a direct consequence. 

Turning to part~$(ii)$, { let $K$ be an arbitrary convex
strongly compact subset of $\operatorname{supp}({\rm P})$}. Starting from $K_0 \coloneqq K$, one can create an increasing sequence $\{K_i\}_{i=0}^{\infty}$ of strongly compact convex\footnote{Without loss of generality~we can assume that $K_i$ is convex for all $i\in \N$; otherwise we replace each $K_i$ by its closed convex hull, which still enjoys strong compactness (see~\citet[Theorem 3.20]{RudinFunction}).}  sets such that ${\rm P}(\H\setminus K_i)\leq \frac{1}{2^i}$. 
Let $\mathbb{K} \coloneqq\bigcup_{i\in \N} K_i$. The functions 
$\psi_n$ can be extended to be continuous on~$\H$ in view of the fact that $\bigcup_{x\in \H}\partial \psi_n(x)\subset \mathcal{B}(0,M)$ (see e.g.,~\eqref{eq:countinous}). Let $\mathcal{C}(\mathbb{K})$ denote the space of bounded continuous functions endowed with the metric
$$ \|f\|_{\mathbb{K}} \coloneqq\sum_{j=0}^{\infty} \frac{\|f\|_{K_j} }{(\operatorname{diam}(K_j)+1)^{2j}}$$
where $\|f\|_{K_j} \coloneqq{\sup_{x\in K_j}  |f(x)|}$ and $\operatorname{diam}(K)\coloneqq \sup_{x,y \in K} \|x-y\|$. Observe that $f_n\to f$ in $\mathcal{C}(\mathbb{K})$ if and only if $f_n\to f$ uniformly in $K_i$ for all $i\in \N$. Since $\psi$, in part $(ii)$ of the theorem,  is unique  up to additive constants only, we can set $\psi_n(x_0)=\psi(x_0)=0$ for some $x_0\in {\rm supp}({\rm P})$.  Note that, for any~$x^1,x^2\in K_0$ and $y_n^1\in \partial \psi_n(x^1)$, $y_n^2\in \partial \psi_n(x^2)$, the inequalities
\begin{equation}
    \label{ineqProofStab}
    \psi_n(x^1)- \psi_n(x^2)\leq \langle y_n^1,x^1-x^2\rangle \quad \text{and}\quad 
 \psi_n(x^1)- \psi_n(x^2) \geq \langle y_n^2, x^1-x^2\rangle,
\end{equation}
along with the assumption of bounded support, imply $ |\psi_n(x^1)- \psi_n(x^2)|\leq M\| x^1-x^2\|$ for~$n$ large enough.
Then, the Arzel\`a-Ascoli theorem implies the uniform convergence of $\psi_n$ (along a subsequence $\{ n^0_k\}_{k\in \N}$) to a function $\rho_0$ in $K_0$. Set $n_0\in \N$ such that $\| \rho_0-{\psi_{n_0}}\|_{K_0}\leq  {1}/{2^{0} =1}$.    We construct a general $\rho$ by using a diagonal argument: for $K_1$, there exists a subse-\linebreak quence $\{ n^1_k\}_{k\in \N}$ of  $\{ n^0_k\}_{k\in \N}$ such that $\psi_{n^1_k}$ converges uniformly in $K_1$ to a function $\rho_1$. Set $n_1\in \N$ such  that~$\| \rho_1-\psi_{n_1}\|_{K_1}\leq {1}/{2^1}=1/2$.  Note that $\rho_1$ agrees with $\rho_0$ in $K_0$. Continuing in this fashion  for~$j=1,2,\ldots $,  define $\rho_j$ on $K_j$, which agrees with $\rho_i$ on $K_i$ for $i < j$ with
\begin{equation}\label{eq:Rho-Psi}
\| \rho_j-\psi_{n_j}\|_{K_j}\leq \frac{1}{2^j}.
\end{equation}
This iterative construction yields  $\rho$ as the unique function in $\mathbb{K}$ which agrees with $\rho_i$ in  $K_i$ for all~$i\in{\mathbb N}$ and   a sequence  $\{n_i\}_{i\in \N}$ such that $$ \|\psi_{n_i}-\rho\|_{\mathbb{K}}\leq \frac{C}{2^i}+\sum_{j=n_i}^{\infty}\frac{\|\rho_j\|_{K_j}\|\psi_{n_i}\|_{K_j}}{(\operatorname{diam}(K_j)+1)^{2j}}$$
where $C \coloneqq\sum_{j=0}^{\infty} {1}/{(\operatorname{diam}(K_j)+1)^{2j}}<\infty$.
Since $$\|\rho_j\|_{K_j}\leq \| \psi_{n_j}\|_{K_j}+ {1}/{2^j}\leq M \operatorname{diam}(K_j) + {1}/{2^j},$$  the rest 
$\sum_{j=n_i}^{\infty} {\|\rho_j\|_{K_j}\|\psi_{n_i}\|_{K_j} }/{(\operatorname{diam}(K_j)+1)^{2j}}$ of the series  tends to zero as $i \to \infty$. As a consequence, 
$ \| \psi_{n_i}-\rho\|_{\mathbb{K}} \to 0$ as $i\to \infty$.  
We claim that this limit $\rho$ is convex and continuous in~$\mathbb{K}$. To prove  convexity,  set $t\in (0,1)$, $(x,y)\in \mathbb{K}^2$, and take limits (as $i \to \infty$) in
$$ \psi_{n_i}(tx+(1-t)y)\leq t\psi_{n_i}(x)+ (1-t)\psi_{n_i}(y). $$
To prove  continuity in $\mathbb{K}$,  set $(x^1, x^2)\in \mathbb{K}^2$. There exists some $K_j$ such that $x^1$ and $x^2$ both lie in~$K_\ell$ for $\ell\geq j$. Since $ |\psi_{n_i}(x^1)- \psi_{n_i}(x^2)|\leq M\| x^1-x^2\|$ and $\|\psi_{n_i}-\rho\|_{K_j}\to 0 $, we conclude that $$ |\rho(x^1)- \rho(x^2)|\leq M\| x^1-x^2\|$$ by letting $i\to \infty$ in $|\rho(x^1)- \rho(x^2)|\leq 2\|\rho- \psi_{n_i}\|_{K_j}+|\psi_{n_i}(x^1)- \psi_{n_i}(x^2)|.$
Set $(i,j)\in \N^2$ and characterize the measure $\gamma_{n_i}^{K_j}$ by imposing $\int f(x,y)d\gamma_{n_i}^{K_j}(x,y) =\int \mathbb{I}_{K_j}(x)f(x,y)d\gamma_{n_i}(x,y)$ (where~$\mathbb{I}_{K_j}$ denotes the indicator function of the set $K_j$)  for any continuous bounded\linebreak  function~$f:\H\times \H\to \R$. Denote by $ {\rm P}_{n_i}^{K_j}$ and $ {\rm Q}_{n_i}^{K_j}$ its marginals. Since $\operatorname{supp}(\gamma_{n_i})$ is cyclically monotone, $\operatorname{supp}(\gamma_{n_i}^{K_j})$ is cyclically monotone as well. The `truncated' conjugate function
$$\psi^*_{n_i,K_j}(y) \coloneqq\sup_{x\in K_j}\{ \langle x,y \rangle -\psi_{n_i}(x)\}, \qquad  \;\; y \in \H$$
satisfies 
\begin{equation}
    \label{gammaA.s.1}
   \psi_{n_i}(x)+\psi^*_{n_i,K_j}(y)=\langle x,y \rangle\quad  \text{for $\gamma_{n_i}^{K_j}$-almost all~$(x,y)$}.
\end{equation}
As $i\to +\infty$ with $j$ fixed,  $\gamma_{n_i}^{K_j}$ tends  weakly to the measure $\gamma^{K_j}$ characterized by $$\int f(x,y)d\gamma^{K_j}(x,y) =\int \mathbb{I}_{K_j}(x)f(x,y)d\gamma(x,y)$$ for any continuous bounded function $f:\H\times \H\to \R$. Let us show that the truncated conjugate of $\rho$, namely,
$$\rho^*_{K_j}(y) \coloneqq\sup_{x\in K_j}\{ \langle x,y \rangle -\rho(x)\}=\sup_{x\in K_j}\{ \langle x,y \rangle -\rho_j(x)\}$$
is the uniform (in $\H $) limit of $\psi^*_{n_i,K_j}$ as $i\to \infty$. To prove this claim,  note that for all  $y\in \H$ and all~$i\geq j$, 
$$ |\psi^*_{n_i,K_j}(y)-\rho^*_{K_j}(y)|\leq \|\rho_j-\psi_{n_i}\|_{K_j}=\|\rho_i-\psi_{n_i}\|_{K_j}\leq \|\rho_i-\psi_{n_i}\|_{K_i}\leq \frac{1}{2^i} \xrightarrow{i\to \infty}0 , $$
where the last inequality holds by construction (see \eqref{eq:Rho-Psi}. 
 As a consequence, from Markov's inequality,
$${\rm Q}^{K_i}_{n_i}( \{ y: \ | \rho_{K_j}^*(y)-\psi^*_{n_i,K_j}(y)|>\delta \})\leq \frac{1}{2^i\delta} \to 0 $$ {as $i\to \infty$} with $j$ fixed, for all~$\delta>0$. 
On the other hand, still using Markov's inequality, we obtain
\begin{align*}
    {\rm P}_{n_i}^{K_j}(\{ x:\ |\psi_{n_i}(x)-\rho(x) |>\delta\})&= {\rm P}_{n_i}( \{ x:\ |\psi_{n_i}(x)-\rho_i(x) |>\delta\}\cap K_j) \leq \frac{1}{2^i\delta}, 
\end{align*}
which tends to zero as $i\to \infty$ with $j$ fixed. Then, 
$$\gamma_{n_i}^{K_j} \big(\{ (x,y): \ | \rho(x)+\rho_{K_j}^*(y)-(\psi_{n_i}(x)+\psi^*_{n_i,K_j}(y))|>\delta \}\big)\to 0  $$
as $i\to \infty$ so that, via \eqref{gammaA.s.1},
$$ \gamma_{n_i}^{K_j}(\{ (x,y): \ | \langle x,y\rangle-(\rho(x)+\rho_{K_j}^*(y))|>\delta \})\to 0  $$
as $i\to \infty$. Since (by the continuity of $\rho$  and $\rho_{K_j}^*$ in  $K_j$ and $\H$, respectively) the set $$\{(x,y)\in K_j\times \H:\ | \langle x,y\rangle-(\rho(x)+\rho_{K_j}^*(y))|>\delta\}$$ is open in $K_j\times \H$,  the Portmanteau theorem yields 
\begin{multline*}
    0=\liminf_{i} \gamma_{n_i}^{K_j}(\{ (x,y): \ | \langle x,y\rangle-(\rho(x)+\rho_{K_j}^*(y))|>\delta \})\\
    \geq \gamma^{K_j}(\{ (x,y): \ | \langle x,y\rangle-(\rho(x)+\rho_{K_j}^*(y))|>\delta \}).
\end{multline*}
From this we  conclude that 
$ \gamma^{K_j}(\{ (x,y): \ | \langle x,y\rangle-(\rho(x)+\rho_{K_j}^*(y))|>\delta \})=0 $.
Thus, on the one hand,  
\begin{equation}
\label{a.s.rhosolves}
    \langle x,y\rangle-(\rho_i(x)+\rho_{K_j}^*(y))=0,\quad \text{for ${\gamma}^{K_j}$-almost all $(x,y)$}
\end{equation}
and, on the other hand, 
\begin{equation}
\label{a.s.psisolves}
    \langle x,y\rangle-(\psi(x)+\psi^*(y))=0,\quad \text{for ${\gamma}^{K_j}$-almost all  $(x,y)$}.
\end{equation}
Theorem~\ref{Theorem:Maccan} applied to the marginals ${\rm P}^{K_j}$ and ${\rm Q}^{K_j}$  of ${\gamma}^{K_j}$ (rescaled  in order to be probability  measures) yields $\nabla \rho (x) = \nabla \psi (x)$ for $ {\rm P}^{K_j}$-almost all~$x$. That is, $\nabla \rho (x) = \nabla \psi (x)$ for $ {\rm P}$-almost all~$x\in K_j$. As a consequence,
\begin{align*}
    {\rm P}(\{ x\in \H:\ \nabla \rho (x) \neq  \nabla \psi (x)\})&\leq   {\rm P}(\{ x\in \H:\ \nabla \rho (x) \neq \nabla \psi (x)\}\cap K_j)+{\rm P}(\H\setminus K_j)
    \\&={\rm P}(\H\setminus K_j)\leq \frac{1}{2^j}.
\end{align*}
Letting $j\to \infty$, we obtain $ {\rm P}(\{ x\in \H:\ \nabla \rho (x) = \nabla \psi (x)\})=1$.  Then, by the assumed uniqueness,  there exists $a\in \R$ such that $\rho=\psi+a$ $\rm P$-a.s. Such an $a$ must be zero due to the fact that~$\rho(x_0)=\lim_{i\to \infty}\psi_{n_i}(x_0)=\psi(x_0)$. 

To conclude, let us show that $\rho(x)=\psi(x)$ for all $x\in K$. Assume that $\rho(x_1)\neq\psi(x_1)$ for some $x_1\in K_0=K\subset \operatorname{supp}({\rm P})$. The continuity of $\rho-\psi$ in $\mathbb{K}$ implies that $\rho(x)\neq\psi(x)$ for all~$x\in \mathcal{B}(x_1,\epsilon)\cap \mathbb{K}$ with $\epsilon>0$ small enough. This, however, cannot be  since ${\rm P}(\mathcal{B}(x_1,\epsilon)\cap \mathbb{K})>0$. Hence, 
$\rho(x)=\psi(x)$ for all $x\in K$, as was to be shown.   
\end{proof}

\section*{Funding}
The research of Alberto Gonz\'{a}lez-Sanz is partially supported by grant PID2021-128314NB-I00 funded by MCIN/AEI/ 10.13039/501100011033/FEDER, UE and by the AI Interdisciplinary Institute ANITI, which is funded by the French ``{\it investing for the Future – PIA3}" program under
the Grant agreement ANR-19-PI3A-0004. The research of Marc Hallin is supported by the Czech Science Foundation grant 
GA\v{C}R22036365, and the research of Bodhisattva Sen is supported by NSF grant DMS-2015376.


\end{document}